\documentclass{amsart}

\usepackage[utf8]{inputenc}
\usepackage{amsmath,amssymb,amsthm,units, stmaryrd,stackrel,relsize,bm}

\newtheorem{theorem}{Theorem}[section]

\newtheorem{lemma}[theorem]{Lemma}
\newtheorem{sublemma}[theorem]{Sublemma}

\newtheorem{corollary}[theorem]{Corollary}
\newtheorem{proposition}[theorem]{Proposition}
\newtheorem{question}[theorem]{Question}

\theoremstyle{definition}
\newtheorem{definition}[theorem]{Definition}

\theoremstyle{remark}
\newtheorem{remark}[theorem]{Remark}

\usepackage[
    colorlinks=true, 
    citecolor=red,
    linkcolor=blue,
    urlcolor=blue]{hyperref}
\newcommand\ad{{\mathsf{AD}}}

\newcommand\zfc{{\mathsf{ZFC}}}
\newcommand\zf{{\mathsf{ZF}}}
\newcommand\kp{{\mathsf{KP}}}
\newcommand\dc{{\mathsf{DC}}}

\newcommand\Ord{{\mathsf{Ord}}}

\newcommand\dom{{\text{dom}}}

\newcommand\cof{{\text{cof}}\,}

\newcommand\coll{\text{Coll}}
\newcommand\rep{\text{Rep}}

\newcommand{\PI}{\boldsymbol\Pi}
\newcommand{\SIGMA}{\boldsymbol\Sigma}
\newcommand{\DELTA}{\boldsymbol\Delta}

\newcommand\quotesleft{``}
\newcommand\quoteleft{`}
\makeatletter
\def\Ddots{\mathinner{\mkern1mu\raise\p@
\vbox{\kern7\p@\hbox{.}}\mkern2mu
\raise4\p@\hbox{.}\mkern2mu\raise7\p@\hbox{.}\mkern1mu}}
\makeatother

\begin{document}
\title{Long Borel Games}
\author{J. P. Aguilera}
\address{Institute of Discrete Mathematics and Geometry, Vienna University of Technology. Wiedner Hauptstra{\ss}e 8--10, 1040 Vienna, Austria.}
\email{aguilera@logic.at}

\begin{abstract}
We study games of length $\omega^2$ with moves in $\mathbb{N}$ and Borel payoff.  These are, e.g., games in which two players alternate turns playing digits to produce a real number in $[0,1]$ infinitely many times, after which the winner is decided in terms of the sequence belonging to a Borel set in the product space $[0,1]^\mathbb{N}$.

The main theorem is that Borel games of length $\omega^2$ are determined if, and only if, for every countable ordinal $\alpha$, there is a fine-structural, countably iterable model of Zermelo set theory with $\alpha$-many iterated powersets above a limit of Woodin cardinals.
\end{abstract}
\maketitle

\setcounter{tocdepth}{1}
\tableofcontents
\section{Introduction}
Given a Borel subset $A$ of the Hilbert space $\ell^2$, define a two-player perfect-information, zero-sum game in which Player I and Player II alternate turns playing digits $x^0_i \in \{0,1,\hdots,9\}$ infinitely often, giving rise to a real number 
$$x^0 = \sum_{i = 0}^\infty \frac{x^0_i}{10^{i+1}}.$$
Afterwards, they continue alternating turns to produce real numbers $x^1, x^2$, and so on. To ensure that the sequence belongs to $\ell^2$,
we normalize, e.g., by setting
$$y^i = \frac{x^i}{i+1}.$$
Player I wins the game if, and only if, $\vec y \in A$; otherwise Player II wins.
These are games of transfinite length $\omega^2$. (Recall that $\omega$ denotes the order-type of the natural numbers.)
It is a natural question whether these games are \emph{determined}, i.e., whether one of the players has a winning strategy. 

Given $A\subset [0,1]$, the Gale-Stewart game on $A$ is defined as follows: two players, I and II, alternate infinitely many turns playing digits $x(0), x(1), \hdots$, as before. This play gives rise to a real number $x \in [0,1]$ as before.
Player I wins the game if, and only if, $x \in A$; otherwise, Player II wins. We say that a set $A\subset\mathbb{R}$ is \emph{determined} if the game associated to $A \cap [0,1]$ is determined.

Our games of length $\omega^2$ can be recast in terms of Gale-Stewart games: given $A \subset [0,1]$, players I and II alternate infinitely many turns playing digits $x^0_0, x^0_1, x^0_2,\hdots$; afterwards, they alternate infinitely many turns playing digits $x^1_0, x^1_1, x^1_2,\hdots$; and so on. Let 
$$\rho:\mathbb{N}\times\mathbb{N} \to \mathbb{N}$$
be a Borel bijection, e.g., 
$$(i,j)\mapsto \frac{1}{2}(i+j)(i+j+1) + j.$$
Put
$$x = \sum_{j = 0}^\infty\sum_{i = 0}^\infty \frac{x^j_i}{10^{\rho(i,j)+1}}.$$
Player I wins the game if, and only if, $x \in A$; otherwise, Player II wins.

A third and most useful way of studying these games is in terms of the Baire space $\mathbb{N}^\mathbb{N}$: given a Borel $A \subset\mathbb{N}^\mathbb{N}$, two players alternate $\omega^2$ many turns playing elements of $\mathbb{N}$ and produce an infinite sequence $\{x_i:i\in\mathbb{N}\}$, with $x_i \in\mathbb{N}^\mathbb{N}$. As before, by means of a Borel bijection, this sequence is identified with an element $x\in\mathbb{N}^\mathbb{N}$. Player I wins the game if, and only if, $x\in A$; otherwise, Player II wins. It is this the way in which we will speak of games of length $\omega^2$ hereafter. Every Polish space is a continuous image of $\mathbb{N}^\mathbb{N}$ and Borel sets have enough closure properties that  these forms of determinacy are all equivalent. 

{\begin{center}
{{---}{---}}
\end{center}
}

Zermelo set theory ($\mathsf{Z}$) was introduced in 1908 by E. Zermelo with the goal of reducing Cantor and Dedekind's theories to a small collection of principles. It, together with A. Fraenkel and T. Skolem's refinements, resulted in Zermelo-Fraenkel set theory ($\zf$). In this article, $\mathsf{Z}$ consists of the axioms of extensionality, separation, power set, pairing, union, foundation, and infinity; and $\zf$ consists of $\mathsf{Z}$, plus the axiom(s) of replacement. $\mathsf{ZC}$ and $\zfc$ are the results of adding the Axiom of Choice to $\mathsf{Z}$ and $\zf$, respectively.

Although $\mathsf{ZC}$ is a powerful theory in which a great deal of mathematics can be formalized, there are many examples of theorems that cannot be proved without the use of the axiom of replacement. An early example of this was that of \emph{Borel determinacy}, i.e., whether all Borel subsets of $\mathbb{R}$ are determined.

\begin{theorem}[Martin \cite{Ma75}]\label{BDT}
Suppose $A$ is Borel. Then $A$ is determined.
\end{theorem}

Even before Martin's Borel determinacy theorem had been proved, it was known that any proof would need essential use of the axiom of replacement:

\begin{theorem}[Friedman \cite{Fr71}]\label{harveystheorem}
Borel determinacy is not provable in $\mathsf{ZC}$.
\end{theorem}

Friedman's work showed that any proof of Borel determinacy would require the use of arbitrarily large countable iterations of the power set operator, and Martin's work showed that this suffices. A convenient slogan is that Borel determinacy captures the strength of countably iterated powersets. 

It is known that one cannot hope to prove the determinacy of any reasonable class of non-Borel games within $\zfc$; however, one can extend $\zfc$ by very natural \emph{large cardinal axioms}, or strengthenings of the axiom of infinity, and use them to prove the determinacy of games on the natural numbers. These connections between large cardinal axioms and infinite games have been the focus of extensive study for decades; some of the most famous results include Martin's \cite{Ma70} proof of analytic determinacy, Martin and Steel's \cite{MaSt89} proof of projective determinacy, and Woodin's \cite{Wo88} proof of the Axiom of Determinacy in $L(\mathbb{R})$, the smallest model of $\zf$ containing all reals and all ordinals. (Recall that the Axiom of Determinacy is the assertion that every Gale-Stewart game (of length $\omega$) is determined. By a result of Mycielski and Steinhaus, it is inconsistent with the Axiom of Choice.)
More specifically, Woodin showed that if there are infinitely many Woodin cardinals and a measurable cardinal greater than them, then the Axiom of Determinacy holds in $L(\mathbb{R})$. The measurable cardinal is necessary, in the sense that the existence of infinitely many Woodin cardinals does not imply that the Axiom of Determinacy holds in $L(\mathbb{R})$.

Our main theorem shows that Borel determinacy for games of length $\omega^2$ captures the strength of countably iterated powersets \emph{above a limit of Woodin cardinals}.\footnote{In the interest of expedience, one may define  the statement \quotesleft there are infinitely many Woodin cardinals'' to mean the following in $\mathsf{ZC}$: there is an infinite sequence $(\delta_0,\delta_1,\hdots)$ of cardinals such that for each $i\in\mathbb{N}$, $V_{\delta_i}$ is a model of $\zfc$ and $\delta_i$ is a Woodin cardinal in $V_{\delta_j}$ whenever $i<j$.}
\begin{theorem}\label{mainborelwoodins}
The following are equivalent over $\zfc$:
\begin{enumerate}
\item Borel determinacy for games of length $\omega^2$;
\item For every countable $\alpha$, there is a countably iterable extender model of $\mathsf{ZC}$ satisfying ``$V_{\lambda+\alpha}$ exists, where $\lambda$ is a limit of Woodin cardinals.''
\end{enumerate}
\end{theorem}
Hence, Borel determinacy for games of length $\omega^2$ is a natural example of a theorem whose strength can only be correctly gauged over $\zfc$ in terms of weaker set theories.
The models referred to in Theorem \ref{mainborelwoodins} are fine-structural and built over some $x\in\mathbb{R}$.
As a consequence thereof, we obtain a result on the provability of Borel determinacy for games of length $\omega^2$ similar to the one for short games:
\begin{corollary}\label{mainborelwoodinscor}
Borel determinacy for games of length $\omega^2$ is provable in the theory $\zfc$ + ``there are infinitely many Woodin cardinals,'' but not in the theory $\mathsf{ZC}$ + ``there are infinitely many Woodin cardinals.''
\end{corollary}
The first part of Corollary \ref{mainborelwoodinscor} follows from Theorem \ref{mainborelwoodins}, together with a theorem of Steel \cite{St93} and a reflection argument. Our proof of Theorem \ref{mainborelwoodins} does not quite go through if one replaces $\zfc$ by $\mathsf{ZC}$ in the statement, although we conjecture that the equivalence is still provable there. Nonetheless, the proof does go through if one replaces $\zfc$ by $\mathsf{ZC}$, together with the assumption that there are enough ordinals\footnote{It suffices to assume that $\vartheta_\alpha$ exists for all $\alpha<\omega_1$ (see below).} and, in particular, if one assumes that there are infinitely many Woodin cardinals. This yields the second part of Corollary \ref{mainborelwoodinscor}.

The main step towards proving Theorem \ref{mainborelwoodins} is a ``shortening'' result like the one in \cite{Ag18a}.
It shows that Borel games of length $\omega^2$ are determined if, and only if, those in a much more complicated class of games of length $\omega$ are. Recall the definition of G\"odel's constructible hierarchy over $\mathbb{R}$: $L_0(\mathbb{R})$ is $V_{\omega+1}$ (the collection of all sets of von Neumann rank ${\leq}\omega$), $L_{\alpha+1}(\mathbb{R})$ is the collection of all sets definable over $L_{\alpha}(\mathbb{R})$ from parameters in $L_{\alpha}(\mathbb{R})$, and $L_\lambda(\mathbb{R}) = \bigcup_{\alpha<\lambda}L_{\alpha}(\mathbb{R})$ if $\lambda$ is a limit ordinal.

\begin{theorem}\label{theorem1intro}
Let $\vartheta$ be the least ordinal such that for every countable ordinal $\alpha$, there is $\xi<\vartheta$ such that $L_\xi(\mathbb{R})$ satisfies \quotesleft $V_\alpha$ exists.''
The following are equivalent:
\begin{enumerate}
\item Borel determinacy for games of length $\omega^2$;
\item $L_\vartheta(\mathbb{R})$ is a model of the Axiom of Determinacy.
\end{enumerate}
\end{theorem}
It is independent of $\zfc$ whether $L_\vartheta(\mathbb{R})$ is a model of the Axiom of Determinacy or of the Axiom of Choice (or neither), but it is never a model of Zermelo set theory, since $\mathcal{P}(\mathbb{R}) \cap L_\vartheta(\mathbb{R})\not\in L_\vartheta(\mathbb{R})$.

The use of large cardinals is essential in proofs of determinacy by G\"odel's second incompleteness theorem, for one can often prove the existence of inner models of set theory satisfying forms of the axiom of infinity from the determinacy of infinite games. 
Results of this form include Friedman's Theorem \ref{harveystheorem} and Theorem \ref{mainborelwoodins}, but also Harrington's \cite{Ha78} work on analytic determinacy and Woodin's work in projective determinacy and the Axiom of Determinacy (see e.g., \cite{MSW} and \cite{KW10}). More recent work in this direction includes that of Steel \cite{St}, Montalb\'an-Shore \cite{MoSh11}, Welch \cite{We11}, Trang \cite{Tr13}, and others. In order to obtain equivalences such as Theorem \ref{mainborelwoodins}, as opposed to equiconsistency results, one needs to involve the techniques of the core model induction. Since we shall be confined within $L(\mathbb{R})$, the arguments we will need can be found in Steel-Woodin \cite{StW16}. Other expositions include Schindler-Steel \cite{SchSt} and Wilson \cite{Wi12}.

The determinacy of games of countable length with analytic (or projective) payoff has been proved by Neeman \cite{Ne04} from assumptions that are very likely optimal. The optimality for games of sufficiently closed length has been verified by Trang \cite{Tr13} and by Woodin, in unpublished work.
Reversals of determinacy hypotheses for projective games of length $\omega^2$ were obtained in \cite{AgMu} and for clopen games in \cite{Ag18a}. Games of length $\omega^2$ were first considered in print by Blass \cite{Bl75}. Let us finish this introduction with some open questions:

\begin{question}
Is the equivalence in Theorem \ref{mainborelwoodins} provable within $\mathsf{ZC}$?
\end{question}

\begin{question}
What is the consistency strength of Borel determinacy for games of length $\omega^2+\omega$?
\end{question}

\begin{question}
What is the consistency strength of determinacy for games of length $\omega$ with payoff in the smallest $\sigma$-algebra containing the projective sets?
\end{question}

{\begin{center}
{{---}{---}}
\end{center}
}

The article is structured as follows: Theorem \ref{theorem1intro} is proved in Section \ref{BDetLR}.
The proof is purely descriptive-set-- and recursion-theoretic. It involves the use of model games which incorporate elements from the constructions in Martin and Steel \cite{MaSt08}, Friedman \cite{Fr71}, and Martin (see \cite{Ma}).
A tool used in the proof (Lemma \ref{TheoremDG}) is a method of describing  initial segments of the $L(\mathbb{R})$-hierarchy by games in which two players determine the truth of a statement in $L(\mathbb{R})$ by challenging each other's claims with auxiliary games. This lemma is proved in Section \ref{SectDG}.
The article ends with Section \ref{BDSectDMT}, in which localizations of arguments from inner model theory are shown to imply Theorem \ref{mainborelwoodins}, starting from Theorem \ref{theorem1intro}.\\

\paragraph{\textbf{Acknowledgements}} We would like to thank Sandra M\"uller, Grigor Sargsyan, and Hugh Woodin for engaging in fruitful discussions with us. This work was partially suported by a grant from the Austrian Science Fund.

\section{Preliminaries}
We shall use notions from descriptive set theory and ordinal recursion theory freely. We refer the reader to Moschovakis' book \cite{Mo09} and to Barwise's book \cite{Ba75} for background. Additional background, as well as history, can be found in the chapters from the Handbook of Set Theory by Schindler-Zeman \cite{SchZe10}, Steel \cite{St10}, and Koellner-Woodin \cite{KW10}.

Our notation is standard.
Among the usual abuses of notation in which we shall engage are the identification of $\mathbb{R}$ with the Baire space $\mathbb{N}^\mathbb{N}$, as well as the identification of Gale-Stewart games with other infinite zero-sum games without explicitly transforming the payoff set into that of a Gale-Stewart game. We also identify $\mathbb{R}$ with $V_{\omega+1}$. In particular, natural numbers are real numbers. 

If $\sigma$ is a strategy for a short game (i.e., one of length $\omega$) on $\mathbb{N}$, then we may alternately regard $\sigma$ as a real, as a function on $\mathbb{N}^{\mathbb{N}}$, or as the function on $\mathbb{R}$ which to each $x\in\mathbb{R}$ assigns the unique $y$ which results by playing against $x$ in a way consistent with $\sigma$. If $\sigma$ and $\tau$ are two strategies, $\sigma*\tau$ denotes the result of facing them off against each other.

We will often need to consider games of length $\omega$ with moves in $\mathbb{R}$. These can be considered as games of length $\omega^2$ on $\mathbb{N}$ in which Player I's moves are ignored throughout the $(2n+1)$th block of $\omega$-many moves, and Player II's moves are ignored throughout the $(2n)$th block of $\omega$-many moves. In particular, if Borel games of length $\omega^2$ are determined, then so too are Borel games of length $\omega$ with moves on $\mathbb{R}$; the converse is not true, as the latter statement is provable in $\zfc$ (by Martin's proof of Borel determinacy). 

At no point other than Section \ref{BDSectDMT} will we need any inner model theory, but we will make use of simple fine-structural facts about $L(\mathbb{R})$. We make use of the Jensen hierarchy over the reals (see \cite{Je72}). Recall its definition given by:
\begin{align*}
J_0(\mathbb{R}) 
&= V_{\omega+1}\\
J_{\alpha+1}(\mathbb{R}) 
&= \text{closure of $J_\alpha(\mathbb{R}) \cup \{J_\alpha(\mathbb{R})\}$ under rudimentary functions} \\
J_\lambda(\mathbb{R}) 
&= \bigcup_{\alpha<\lambda}J_\alpha(\mathbb{R}), \text{ for limit $\lambda$}.
\end{align*}

The following definition is perhaps not standard, but will suffice for our purposes:
\begin{definition}
Let $\alpha$ and $\beta$ be ordinals. We say that the \emph{projectum} of $J_\alpha(\mathbb{R})$ is $V_\beta$ if $\beta<\alpha$ and $\beta$ is the least ordinal such that
\[V_{\beta+1}\cap \big(J_{\alpha+1}(\mathbb{R}) \setminus J_\alpha(\mathbb{R})\big)\neq\varnothing.\]
\end{definition}
Let $\beta$ be such that $V_\beta$ is the projectum of $J_\alpha(\mathbb{R})$. We say that $J_\alpha(\mathbb{R})$ \emph{projects to} $V_\beta$. Since $\mathbb{R}$ always belongs to $J_\alpha(\mathbb{R})$, $\beta$ is at least $\omega+1$. If $\beta = \omega+1$, we say that $J_\alpha(\mathbb{R})$ projects to $\mathbb{R}$. If $J_\alpha(\mathbb{R})$ does not project to $\mathbb{R}$, then every subset of $\mathbb{R}$ definable over $\mathcal{P}(\mathbb{R})^{J_\alpha(\mathbb{R})}$ belongs to $J_\alpha(\mathbb{R})$, so
\[\big(\mathbb{N},\mathbb{R},\mathcal{P}(\mathbb{R}) \cap J_\alpha(\mathbb{R})\big)\]
is a model of third-order arithmetic.

If the projectum of $J_\alpha(\mathbb{R})$ is not defined, then $J_\alpha(\mathbb{R})$ contains all definable subsets of $V_\beta^{J_\alpha(\mathbb{R})}$ for all $\beta<\alpha$, and so 
\[J_\alpha(\mathbb{R})\models\zf.\]

\begin{proposition}\label{PropositionProjectum}
Suppose that the projectum of $J_\alpha(\mathbb{R})$ is $V_{\beta}$ and $\omega+1<\beta$. Then, there is a surjection
\[\rho:V_{\beta}\cap J_{\alpha}(\mathbb{R})\twoheadrightarrow J_\alpha(\mathbb{R})\]
in $J_{\alpha+1}(\mathbb{R})$.
\end{proposition}
\proof
This proof is similar to those in \cite[Section 3]{Je72} and \cite[Section 1]{St08}. 
Suppose that the projectum of $L_\alpha(\mathbb{R})$ is $V_\beta$, with $\omega+1<\beta$, and let $\phi$ be a formula in the language of set theory such that for some $a \in J_\alpha(\mathbb{R})$, we have
\[\big\{x \in J_\alpha(\mathbb{R}): J_\alpha(\mathbb{R})\models\phi(a,x)\big\}\in V_{\beta+1}\cap \big(J_{\alpha+1}(\mathbb{R}) \setminus J_\alpha(\mathbb{R})\big).\]
Call the set just defined $A$. 
Suppose that $\phi$ is $\Sigma_n$. In $L(\mathbb{R})$, every set is definable from a real and ordinal parameters, in fact, by Lemma 1.4 of \cite{St08}, there are surjections 
\[f_\beta: [\omega\cdot\beta]^{<\omega}\times\mathbb{R}\twoheadrightarrow J_\beta(\mathbb{R})\]
which are uniformly $\Sigma_1$ over $J_\beta(\mathbb{R})$. Thus, we may assume that $a \in\mathbb{R}$ by replacing $n$ with a bigger natural number if necessary, and replacing $\phi(a,x)$ with the formula
\begin{align*}
\phi(f_\alpha(\vec \beta,a),x)\text{ holds for the least tuple $\vec \beta$ such that the set }\\
\big\{x: \phi(f_\alpha(\vec \beta,a),x)\big\} \text{ does not exist}. \qquad
\end{align*}

Since $\omega+1 < \beta$, Lemma 1.7 and Theorem 1.16 of \cite{St08} imply that $J_\alpha(\mathbb{R})$ has a $\SIGMA_n$ Skolem function, i.e., that there is a partial map
\[h: J_\alpha(\mathbb{R})\times\mathbb{R}\to J_\alpha(\mathbb{R})\]
which is $\Sigma_n$ in a parameter $p \in J_\alpha(\mathbb{R})$ such that whenever $S$ is a nonempty $\Sigma_n$ set over $J_\alpha(\mathbb{R})$ with a parameter $q \in J_\alpha(\mathbb{R})$, then
\[h(q,x) \in S \text{ for some } x\in\mathbb{R}.\]
Thus, one can define the $\Sigma_n$ Skolem hull of $J_\alpha(\mathbb{R})$ with a $\Sigma_n$ formula using parameters from $J_\alpha(\mathbb{R})$. Let $p$ be the parameter defining $h$ and
\[\mathcal{H} = \text{Hull}^{J_{\alpha}(\mathbb{R})}_{\Sigma_{n}}(V_\beta\cap J_\alpha(\mathbb{R}), \{p\}).\]
By Lemma 1.9 of \cite{St08},
\begin{enumerate}
\item $\mathcal{H}$ is a $\Sigma_{n}$-elementary substructure of $J_\alpha(\mathbb{R})$, and
\item $\mathcal{H} = f[X^{<\omega}]$ for some partial map $f$ which is $\Sigma_n$ over $\mathcal{H}$ with parameters in $\mathcal{H}$, where $X = (V_\beta\cap J_\alpha(\mathbb{R})) \cup \{p\}$.
\end{enumerate}
The latter conclusion easily implies that in fact 
\[\mathcal{H} = f[X^{<\omega}]\] 
for some partial map $f$ as above, where $X = V_\beta\cap J_\alpha(\mathbb{R})$.

Since $\omega+1<\beta$, $\mathcal{H}$ contains all reals, so, by condensation, there is $\gamma\leq\alpha$ such that the transitive collapse of $\mathcal{H}$ is $J_\gamma(\mathbb{R})$. Let
\[\pi:J_\gamma(\mathbb{R}) \twoheadrightarrow \mathcal{H}\]
be the collapse embedding, which has critical point ${\geq}\beta$. By definition, $A \in V_{\beta+1}$, so
\begin{align*}
\text{if $J_\alpha(\mathbb{R})\models \phi(a,x)$, then $x \in V_\beta$.}
\end{align*}
Since $V_\beta\cap J_\alpha(\mathbb{R})\subset J_\gamma(\mathbb{R})$ and $a \in \mathbb{R}$, we have for every $x\in J_\gamma(\mathbb{R})$,
\begin{align*}
J_\alpha(\mathbb{R})\models \phi(a,x) 
&\text{ if, and only if, } J_\alpha(\mathbb{R})\models x \in V_{\pi(\beta)}\wedge \phi(a,x),\\ 
&\text{ if, and only if, } J_\gamma(\mathbb{R})\models x \in V_{\beta}\wedge \phi(a,x).
\end{align*}
This is not entirely straightforward, as $V_\beta\cap J_\alpha(\mathbb{R})$ need not be an element of $J_\gamma(\mathbb{R})$ (or of $J_\alpha(\mathbb{R})$), but for each $\eta$, the set
\[\big\{x \in J_\eta(\mathbb{R}): x \text{ has rank } \beta\big\}\]
is uniformly $\Delta^{J_\eta(\mathbb{R})}_1(\beta)$ (provided $\beta<\eta$). Thus, the set $A$ is definable over $J_\gamma(\mathbb{R})$, which implies that $\gamma = \alpha$.  Using $\pi$ and the function $f$ from above, one obtains a partial surjection from $V_\beta\cap J_\alpha(\mathbb{R})$ onto $J_\alpha(\mathbb{R})$ which is definable over $J_\alpha(\mathbb{R})$ by a $\Sigma_n$ formula with parameters. This easily translates into a $\Pi_{n+1}$ (with parameters) total surjection $f^*$ from $V_\beta\cap J_\alpha(\mathbb{R})$ onto $J_\alpha(\mathbb{R})$. Clearly $f^*$ belongs to $J_{\alpha+1}(\mathbb{R})$, which completes the proof.
\endproof

\section{Determinacy and $L(\mathbb{R})$}\label{BDetLR}
If $A\subset\mathbb{R}\times\mathbb{R}$ is a set and $x\in\mathbb{R}$, we define
$A_x = \{y\in\mathbb{R}: (x,y) \in A\}$
and
\begin{align*}
\Game^\mathbb{R}A 
&= \big\{x\in\mathbb{R}: \text{ Player I has a winning strategy in the game }\\
&\quad \qquad  \text{ with payoff $A_x$ and moves in $\mathbb{R}$}\big\}.
\end{align*}
The set $\Game^\mathbb{R}A$ is defined similarly for $A$ a subset of other spaces, such as $\mathbb{R}^n$. Letting
$\Game^\mathbb{R}\DELTA^1_1 = \{\Game^\mathbb{R}A: \text{ $A$ is Borel}\}$,
we write
\begin{align*}
\gamma^1_1 = \sup\{\alpha: \text{ there is a prewellordering of $\mathbb{R}$ in $\Game^\mathbb{R}\DELTA^1_1$ of rank $\alpha$}\}.
\end{align*}
The ordinal $\gamma^1_1$ is usually denoted by $\delta_{\Game^\mathbb{R}\DELTA^1_1}$. \index{$\Game^\mathbb{R}\DELTA^1_1$} \index{$\gamma^1_1$}\footnote{Note that $\Game^\mathbb{R}\Delta^1_1$ is a much smaller class than $\Delta_1^{L(\mathbb{R})} = \Delta_{\Game^\mathbb{R}\Pi^1_1}$, the class of all sets $A$ such that both $A$ and its complement belong to $\Game^\mathbb{R}\Pi^1_1$.}

\begin{lemma} \label{TheoremDG}
$\mathcal{P}(\mathbb{R})\cap L_{\gamma^1_1}(\mathbb{R})\subset\Game^\mathbb{R}\DELTA^1_1.$
\end{lemma}
The proof of Lemma \ref{TheoremDG} will be delegated to the following section. For now, let us assume that it is true.

\begin{remark}
Under determinacy, Lemma \ref{TheoremDG} can be proved easily by appealing to general Wadge theory. The proof that we shall provide, however, is direct, and goes through in e.g., $\zf + \dc$.
\end{remark}

\begin{definition}\index{$\vartheta_\alpha$}
We denote by $\vartheta_\alpha$ the least ordinal $\xi$ such that
$L_\xi(\mathbb{R})$ is a model of 
$$\kp + \text{Separation} + \quotesleft V_\alpha \text{ exists.''}$$
We also define 
$$\vartheta = \sup_{\alpha<\omega_1}\vartheta_\alpha.$$
\end{definition}

The main result of this section is the following: \\
\begin{theorem}\label{mainborel}
The following are equivalent:
\begin{enumerate}
\item \label{mainborel11} Borel determinacy for games of length $\omega^2$;
\item \label{mainborel22} $L_\vartheta(\mathbb{R}) \models \ad$.
\end{enumerate}
\end{theorem}

We need some preliminaries. We will work with structures in the language of set theory with additional constants $\{\dot x_i:i\in\mathbb{N}\}$ (which will be interpreted as real numbers). We fix:
\begin{enumerate}
\item a formula $\theta(\cdot,\cdot,\cdot)$ in the language of set theory defining in $L(\mathbb{R})$ an $\mathbb{R}$-parametrized sequence of wellorderings the union of whose domain is $L(\mathbb{R})$ (e.g., the one given by Lemma 1.4 of \cite{St08}), and
\item injective functions $n(\cdot)$ and $m(\cdot)$ assigning natural numbers to formulae in the extended language in such a way that their ranges are recursive and disjoint and whenever $\dot x_i$ occurs in a formula $\phi$, then $i < m(\phi)$ and $i<n(\phi)$.
\end{enumerate}

Let us begin with an observation:

\begin{lemma}
Suppose that $M$ is a model of $\kp$ and $M$ contains an ordinal $a$ isomorphic to $\alpha$. Then $(V_a)^M$ (which need not belong to $M$) is wellfounded.
\end{lemma}
\proof
Let us identify the wellfounded part of $M$ with its transitive collapse. In particular, we write $\alpha$ for the ordinal $a$ in the statement of the lemma. Although $(V_\alpha)^M$ need not belong to $M$, it is $\Delta_1$-definable (with parameters) over $M$, so given $x \in (V_\alpha)^M$, $M$ can determine the rank of $x$. Thus, if $x \in^M y$ and $y \in^M (V_\alpha)^M$, then
$$\text{rank}^M(x) < \text{rank}^M(y) < \alpha.$$
Since $\alpha$ is wellfounded, $(V_\alpha)^M$ must be too.
\endproof

\begin{lemma}\label{LemmaTalphaisBorel}
For every $\alpha<\omega_1$, let $T_\alpha$ be the set of all (real numbers coding) complete and consistent theories $T$ in the language of set theory with additional constants $\{\dot x_i:i\in\mathbb{N}\}$ such that 
\begin{enumerate}
\item $T$ extends $\kp + \text{\quotesleft $\mathbb{R}$ exists''} + V = L(\mathbb{R}) + \{\dot x_i \in \mathbb{R}:i\in\mathbb{N}\}$,
\item $T$ contains Skolem schemata of the forms
\begin{align*}
\exists x\in\mathbb{R}\, \phi(x) 
&\to \phi(\dot x_{n(\phi)}),\\
\exists x\, \phi(x) 
&\to \exists x\,\exists \beta\in\Ord\,(\theta(\beta,\dot x_{m(\phi)}, x) \wedge \phi(x));
\end{align*}
\item all models of $T$ contain an ordinal isomorphic to $\alpha$.
\end{enumerate}
Then $T_\alpha$ is a Borel subset of $\mathbb{R}$.
\end{lemma}
\proof
Given a complete and consistent theory $T$ in the language of set theory extending $\kp$, we say that $M$ is a \emph{term model} of $T$ if $M$ is the natural $=$-quotient of the collection of all formulae $\phi$ with one free variable such that
\begin{align*}
\text{\quotesleft there is a unique $x$ satisfying $\phi(x)$''} \in T,
\end{align*}
with the membership predicate and the constants $\dot x_i$ having the obvious interpretations. Suppose that $T$ is as in the statement of the lemma and let $M$ be the term model of $T$. It follows from the Tarski-Vaught criterion and the fact that $M$ satisfies the Skolem schemata that $M$ is in fact a model of $T$. Moreover, it embeds into every model of $T$, so every model of $T$ contains an isomorphic copy of $\alpha+1$ if, and only if, $M$ does.

Fix some wellordering $\preceq$ of $\mathbb{N}$ of length $\alpha+1$. Thus, the following are equivalent:
\begin{enumerate}
\item $T \in T_\alpha$;
\item there is a term model $M$ of $T$, and a function $f$ such that
\begin{enumerate}
\item the range of $f$ is the field of $\preceq$,
\item the domain of $f$ is the set of $\in^M$-predecessors of an ordinal of $M$,
\item for all $y\in \dom(f)$, $f(y)$ is the $\preceq$-least number greater than $f(x)$ for every $x \in^M y$;
\end{enumerate}
\item for all term models $M$ of $T$, there is an $M$-ordinal $m \in M$ such that for all functions $f$, if
\begin{enumerate}
\item the range of $f$ is contained in the field of $\preceq$,
\item the domain of $f$ is the set of $\in^M$-predecessors of $(m+1)^M$,
\item for all $y\in \dom(f)$, $f(y)$ is the $\preceq$-least number greater than $f(x)$ for every $x \in^M y$,
\end{enumerate}
then $f(m)$ is the $\preceq$-largest number.
\end{enumerate}

Therefore $T_\alpha$ is Borel.
\endproof

\begin{lemma}\label{LemmaBorelVartheta}
$\vartheta$ is least such that
$$L_\vartheta(\mathbb{R})\models\text{\quotesleft Borel games on $\mathbb{R}$ are determined.''}$$
\end{lemma}
\proof
Clearly $\vartheta$ has this property, for $L_\vartheta(\mathbb{R})$ has initial segments with enough iterated powersets of the real numbers for Martin's proof \cite{Ma75} to go through; moreover, the existence of winning strategies for games on $\mathbb{R}$ is upwards absolute between transitive sets that contain $\mathbb{R}$. 

Let us verify that $\vartheta$ is the least such ordinal. We need to employ model games such as the one of \cite{MaSt08}.
A similar application of the techniques of \cite{MaSt08} was used in \cite{AgMu} to produce iterable models of $\ad$ with Woodin cardinals. Fortunately, our situation is simpler inasmuch as we do not need to construct models with large cardinals (yet); the drawback is that we cannot even ask for wellfoundedness in the payoff. One deals with this issue as in Martin's reversal of Borel determinacy (unpublished, but see the exercise in pp.~53--54 of \cite{Ma} for the case of $\Pi^0_4$ games on $\mathbb{N}$).

Let $\alpha<\omega_1$ be an infinite successor ordinal. We consider the following game:
\begin{align*}
\begin{array}{c|ccccc}
I & x_0, u_0 & & x_2, u_1 & &\hdots\\
II & & x_1, v_0 & & x_3, v_1 & \hdots
\end{array}
\end{align*}
Here, Players I and II play a countable sequence of reals $\{x_i\}_{i\in\mathbb{N}}$, as well as two countable sequences $u,v \in 2^{\mathbb{N}}$. Fix an enumeration $\{\phi_i:i\in\mathbb{N}\}$ of all formulae in the language of set theory with added constants $\{\dot x_i:\in\mathbb{N}\}$. 
We ask of the players that:
\begin{enumerate}
\item $T_I := \{\phi_i: u_i = 1\}$ is a complete, consistent theory in the expanded language, satisfying the following conditions:
\begin{enumerate}
\item $T_I$ extends the theory $\kp$ + Separation $ +$ \quotesleft $\mathbb{R}$ exists'' + $\dc $ + $V = L(\mathbb{R})$;
\item $T_I$ contains the statement $\dot x_i \in\mathbb{R}$ for each $i\in\mathbb{N}$, and it contains the statement $\dot x_i(n) = m$ if, and only if, $x_i(n) = m$;
\item $T_I$ contains the Skolem schemata
\begin{align*}
\exists x\in\mathbb{R}\, \phi(x) 
&\to \phi(\dot x_{n(\phi)}),\\
\exists x\, \phi(x) 
&\to \exists x\,\exists \beta\in\Ord\,(\theta(\beta,\dot x_{m(\phi)}, x) \wedge \phi(x));
\end{align*}
\item We ask that all models of $T_I$ contain an ordinal isomorphic to $\alpha$. If so, then by the proof of Lemma \ref{LemmaTalphaisBorel}, there must be a formula $\phi_\alpha$ in the expanded language that defines the $\alpha$th ordinal in every model of $T_I$. We demand then that $T_I$ contain the axiom 
\begin{align*}\label{eqAxiomZalpha}
\text{\quotesleft if $\phi_\alpha(x)$ holds, then $V_x$ exists,''}
\end{align*} 
as well as 
$$\text{$\phi_\alpha(x) \wedge \phi(a)\to L_a(\mathbb{R})\not\models \kp$ + Separation + \quotesleft $V_x$ \text{ exists,''}}$$
whenever $\phi$ is a formula and the formula \quotesleft the unique $a$ satisfying $\phi(a)$ is an ordinal smaller than $\alpha$'' belongs to $T_I$.
\end{enumerate}
\item $T_{II} := \{\phi_i: v_i = 1\}$ is a complete, consistent theory in the expanded language with the same requirements as above.
\end{enumerate}
If $T_I$ satisfies the properties above and $\mathcal{M}$ is a model of $T_I$, then the definable closure of $\{x_i^\mathcal{M}:i\in\mathbb{N}\}$ in $\mathcal{M}$ is an elementary substructure of $\mathcal{M}$ isomorphic to the term model of $T_I$. Denote this model by $N_I$. In particular, $N_I$ satisfies $\kp$ + Separation + ``$V_\alpha$ exists'' and no proper initial segment of $N_I$ does. Similarly for Player II. If the above conditions are met, Player I wins if, and only if, one of the following holds:
\begin{enumerate}
\item \label{LemmaBorelVartheta2} $N_I$ is isomorphic to an initial segment of $N_{II}$; \emph{or}
\item \label{LemmaBorelVartheta3} there is an ordinal $a$ of $N_I$ such that $L_a(\mathbb{R})^{N_I}$ is isomorphic to an initial segment of $N_{II}$, but $L_{a+1}(\mathbb{R})^{N_I}$ is not.
\end{enumerate}

In either clause, we do not demand that the initial segment of $N_{II}$ be an element of $N_{II}$.

\begin{sublemma}
The winning condition is Borel.
\end{sublemma}
\proof
Let us momentarily reason in  
$$\kp + \text{Separation} + \text{\quotesleft $\mathbb{R}$ exists''} + \dc + V = L(\mathbb{R}).$$ 
Suppose that $J_\xi(\mathbb{R})\not\models$\quotesleft $V_\alpha$ exists,'' for every ordinal $\xi$. Thus, for every $\xi$, the projectum of $J_\xi(\mathbb{R})$ is strictly smaller than $V_\alpha$. By Proposition \ref{PropositionProjectum}, whenever $\xi$ is an ordinal, there is some 
$$a \in V_\alpha\cap J_{\xi+1}(\mathbb{R})$$
such that $a$
codes $J_\xi(\mathbb{R})$. Now, let $N_I$ and $N_{II}$ be as in the definition of the game (i.e., suppose that $T_I$ and $T_{II}$ satisfy the conditions necessary in order for $N_I$ and $N_{II}$ to be defined). Now, $N_I$ and $N_{II}$ might have initial segments that satisfy \quotesleft $V_\alpha$ exists,'' but they certainly do not have any initial segments satisfying \quotesleft $V_{\alpha+1}$ exists,''
so for every initial segment $N$ of either of $N_I$ or $N_{II}$ of (internal) limit length greater than $\alpha$, we have
\begin{equation*}
\text{ for all $\xi\in\Ord^N$, there is } a\in (V_{\alpha+1})^{N} \text{ coding } (J_{\xi}(\mathbb{R}))^{N}.
\end{equation*}
(One could easily get around this issue by removing the condition that the theories satisfy $\kp$ + Separation, but we have adpoted the rule because this is the theory used in the definition of $\vartheta_\alpha$.)

Back in the real world, we verify that conditions \eqref{LemmaBorelVartheta2} and \eqref{LemmaBorelVartheta3} are Borel. We will show that it is Borel to check whether the elements of $(V_\alpha)^{N_I}$ are (up to isomorphism) those of $(V_\alpha)^{N}$, for some initial segment $N$ of $(V_\alpha)^{N_{II}}$. This suffices for \eqref{LemmaBorelVartheta2}, by the above remark and by replacing $\alpha$ with $\alpha+1$.
As in the proof of Proposition \ref{PropositionProjectum}, one should bear in mind that $(V_\alpha)^N$ need not be an element of $N$, but it is definable over $N$ in a Borel way by defining an initial segment of the rank function along a Borel presentation of $\alpha$; viz., a set $x\in N$ has rank $\alpha$ in $N$ if, and only if, $f(x) = \alpha$ for some (all) function(s) $f$ such that
\begin{enumerate}
\item the range of $f$ is $\alpha+1$,
\item the domain of $f$ is a subset of $N$ containing $x$ and is transitive (so, in particular, $\in^N$ is wellfounded in $\dom(f)$),
\item for all $y\in \dom(f)$, $f(y) = \sup\{f(z):z \in^N y\}$. 
\end{enumerate}

Let $f_I$ (in $V$) be a function assigning to each set in $(V_\alpha)^{N_I}$ its rank in $N_I$ and $f_{II}$ be defined similarly. Consider the following Ehrenfeucht-Fra\"iss\'e--like game:
\begin{enumerate}
\item Player II begins by selecting an ordinal $a \in N_{II}$. Define $N = (J_a(\mathbb{R}))^{N_{II}}$. For the rest of the game, the players attempt to determine whether $(V_\alpha)^{N_{I}}$ is equal to $(V_\alpha)^{N}$ up to isomorphism.
\item During turn $1$, Player I selects some $x_0 \in (V_\alpha)^{N_I}$ or some $y_0 \in (V_\alpha)^{N}$. In the first case, let $\xi_0$ be least such that $x_0 \in (V_{\xi_0})^{N_I}$; Player II needs to respond with some $y_0 \in (V_{\xi_0})^{N}$. In the second case, let $\xi_0$ be least such that $y_0 \in (V_{\xi_0})^{N}$; Player II needs to respond with some $x_0 \in (V_{\xi_0})^{N_I}$.
\item During turn $n+1$, assuming $x_n$, $y_n$, and $\xi_n$ have been defined, Player I plays some $\xi_{n+1} < \xi_n$ and either some $x_{n+1} \in (V_{\xi_{n+1}})^{N_I}$ such that $x_{n+1} \in^{N_I} x_n$ or some $y_{n+1} \in (V_{\xi_{n+1}})^{N}$ such that $y_{n+1} \in^{N_{II}} y_n$. In the first case, Player II must respond with some $y_{n+1} \in (V_{\xi_{n+1}})^{N}$ such that $y_{n+1} \in^{N_{II}} y_n$; in the second case, Player II must respond with some $x_{n+1} \in (V_{\xi_{n+1}})^{N_I}$ such that $x_{n+1} \in^{N_I} x_n$.
\item The first player who cannot make a legal move loses the game.
\end{enumerate}
This game can easily be coded by a clopen game of length $\omega$ on $\mathbb{N}$ with parameters $\alpha$, $N_I,N_{II},f_I$, and $f_{II}$; let us denote it by $EF_0(N_I,N_{II},f_I,f_{II})$ to emphasize this. If there is an ordinal $a \in N_{II}$ such that $(V_\alpha)^{N_{I}} = (V_\alpha)^{J_a(\mathbb{R})^{N_{II}}}$, then clearly Player II has a winning strategy. Moreover, if there is no such ordinal, then one can easily use the fact that $(V_\alpha)^{N_I}$ and $(V_\alpha)^{N_{II}}$ are wellfounded to construct a winning strategy for Player I. 

There is a problem with this game, viz., that one may also want to consider initial segments of $N_{II}$ not of the form $J_a(\mathbb{R})$ for $a$ an ordinal of $N_{II}$. For this, we define a game $EF(N_I,N_{II},f_I,f_{II})$ which is like $EF_0(N_I,N_{II},f_I,f_{II})$, except that Player I begins by selecting an ordinal $a_0$ of $N_I$, after which Player II selects $a \in N_{II}$ and they continue as before to determine whether 
$$(V_\alpha)^{J_{a_0}(\mathbb{R})^{N_{I}}} = (V_\alpha)^{J_a(\mathbb{R})^{N_{II}}}.$$ An argument like the preceding one shows that Player II has a winning strategy in $EF(N_I,N_{II},f_I,f_{II})$ if, and only if, $(V_\alpha)^{N_{I}} = (V_\alpha)^{N}$ for some initial segment $N$ of $N_{II}$.

Now, one sees that the following are equivalent:
\begin{enumerate}
\item the elements of $(V_\alpha)^{N_I}$ are those of $(V_\alpha)^{N}$, for some initial segment $N$ of $(V_\alpha)^{N_{II}}$,
\item there is sequence of surjective functions $\{f_\beta:\beta\leq\alpha\}$  and an initial segment $N$ of $N_{II}$ such that, identifying the wellfounded parts of $N_I$ and $N$ with their transitive collapses,
\begin{enumerate}
\item $f_0$ maps the empty set to the empty set,
\item for all $\beta\leq\alpha$, $f_{\beta+1}: (V_\beta)^{N_I}\to (V_\beta)^N$ maps each set $b$ to the set $\{f_\beta(a): a \in^{N_I} b\}$, in the sense of $N$,
\item for all limit $\beta\leq\alpha$, $f_\beta = \bigcup_{\gamma<\beta}f_\gamma$
\end{enumerate}
\item for every pair of functions $f_I$ and $f_{II}$ assigning respectively to each set in $(V_\alpha)^{N_I}$ and  $(V_\alpha)^{N_{II}}$ its rank in $N_I$ and $N_{II}$, Player I does not have a winning strategy in the game $EF(N_I,N_{II},f_I,f_{II})$.
\end{enumerate}
Since $EF(N_I,N_{II},f_I,f_{II})$ is a clopen game, Player I not having a winning strategy is a $\PI^1_1$ condition. It follows that condition \eqref{LemmaBorelVartheta2} is Borel.

Similarly, to prove that \eqref{LemmaBorelVartheta3} is Borel, it suffices to show that it is Borel to check whether there is $a \in N_I$ such that 
\begin{enumerate}
\item $(V_\alpha)^{J_a(\mathbb{R})^{N_I}} \subset (V_\alpha)^{N_{II}}$
\item $(V_\alpha)^{J_{a+1}(\mathbb{R})^{N_I}} \not\subset (V_\alpha)^{N_{II}}$
\end{enumerate}
An argument as above shows that the conjunction of these two statements is Borel. This proves the sublemma.
\endproof

We have shown that the payoff set of the game is Borel. Moreover, it is easily won by Player I in $V$, for she can play the theory of $J_{\vartheta_\alpha}(\mathbb{R}) = L_{\vartheta_\alpha}(\mathbb{R})$ and, e.g., only recursive reals. Recall that ${\vartheta_\alpha}$ is least such that
$$J_{\vartheta_\alpha}(\mathbb{R})\models \text{$\kp$ + \text{Separation} + \quotesleft $V_\alpha$ exists.''}$$
Moreover, as long as Player II plays by copying the theory chosen by Player I (and plays arbitrary reals), then every strategy that ensures that Player I will win demands that she play the theory of $J_{\vartheta_\alpha}(\mathbb{R})$. Thus, the theory of $J_{\vartheta_\alpha}(\mathbb{R})$ can be easily computed from $\mathbb{R}$ and any winning strategy for Player I. 

Now let $\xi$ be an ordinal such that the game is determined in $J_\xi(\mathbb{R})$. Then $J_\xi(\mathbb{R})$ contains a winning strategy $\tau$ for some player. Since $J_\xi(\mathbb{R})$ contains all plays of the game, $\tau$ is also a winning strategy in $V$, so it must be a winning strategy for Player I. Hence, we must have ${\vartheta_\alpha}<\xi$, since $\tau$ cannot belong to $J_{\vartheta_\alpha}(\mathbb{R})$.
Therefore, if $J_\xi(\mathbb{R})$ contains strategies for the game above for each $\alpha<\omega_1$, then it contains every ordinal below $\vartheta$, which is what was to be shown.
\endproof

\begin{lemma}\label{LemmaBorelPWO}
Let $\xi<\vartheta$. Then there is a prewellordering of $\mathbb{R}$ of length greater than $\xi$ definable by a set in $\Game^\mathbb{R}\DELTA^1_1$.
\end{lemma}
\proof
Choose $\alpha<\omega_1$ such that $\xi<\vartheta_{\alpha}$ and let $(i,x), (j,y)$ be two pairs in $\mathbb{N}\times\mathbb{R}$. The ordinal $\vartheta_{\alpha}$ will play no role in the game other than bounding the complexity of the payoff set.
We consider a two-player game $G(\alpha, i,x,j,y)$ given by
\begin{align*}
\begin{array}{c|ccccccccc}
I & x_0, u_0 & &         & &x_2, u_1 & &&&\hdots\\
II &         & &x_1, v_0 &  & &       & x_3, v_1 & &\hdots\\\hline
I &         & & &y_1, t_0  &  &      &  & y_3, t_1 & \hdots\\
II  & &y_0, s_0 &  &        &  &y_2, s_1 & &&\hdots
\end{array}
\end{align*}
Here, I and II play two runs of the game from Lemma \ref{LemmaBorelVartheta} in parallel, except that their roles are reversed in the second game. As before, we suppose fixed an enumeration $\{\phi_i:i\in\mathbb{N}\}$ of all formulae in the language of set theory with added constants $\{\dot x_i:i\in\mathbb{N}\}$.  The rules are:
\begin{enumerate}
\item \label{LemmaBorelPWORule1} If $\phi_j$ is not a formula with one free variable and with no constant symbols other than $\dot x_0$, then Player II loses. 
\item \label{LemmaBorelPWORule2} If $\phi_i$ is not a formula with one free variable and with no constant symbols other than $\dot x_0$, then Player I loses. 
\item \label{LemmaBorelPWORule3} Both players have to obey the rules in each subgame.
\item \label{LemmaBorelPWORule4} If so, then Player II has to win the second subgame. Moreover, she must set $y_0 = y$ and the theory played must contain the formula \quotesleft there is a unique ordinal satisfying $\phi_j$.''
\item \label{LemmaBorelPWORule5} If so, then Player I has to win the first subgame. Moreover, she must set $x_0 = x$ and the theory played must contain the formula \quotesleft there is a unique ordinal satisfying $\phi_i$.''
\end{enumerate}
If all rules thus far have been obeyed, then we declare the winner according to the following rules (here, we identify the theories played with their term models):
\begin{enumerate}
\item \label{LemmaBorelPWOWinCon2} Every real played in the first subgame must be played in the second subgame by Player II; otherwise, Player I wins.
\item \label{LemmaBorelPWOWinCon3} Every real played in the second subgame must be played in the first subgame by Player I; otherwise, Player II wins.
\item \label{LemmaBorelPWOWinCon4} If no winner has been declared thus far, let $N_I$ be the winning model from the first subgame and $N_{II}$ be the winning model from the second subgame (both with $\alpha$ in their wellfounded parts). Then,
\begin{enumerate}
\item We demand that $N_I$ be isomorphic to an initial segment of $N_{II}$; otherwise, Player II wins. 
\item Let $a_I$ be the unique ordinal in $N_I$ satisfying $\phi_i$ and let $a_{II}$ be the unique ordinal in $N_{II}$ satisfying $\phi_j$. Then,
Player I wins if for every (some) isomorphism $j$ from $N_I$ to an initial segment $N$ of $N_{II}$, we have 
$$N_{II}\models j(a_I) < a_{II}.$$ 
If not, then Player II wins.
\end{enumerate}
\end{enumerate}
Let $\vartheta_\alpha^*$ be the order-type of the set of ordinals definable in $L_{\vartheta_\alpha}(\mathbb{R})$ from a real. Given pairs $(i,x), (j,y)$ as above, if there is not a unique ordinal $\eta<\vartheta_{\alpha}$ such that 
$$L_{\vartheta_{\alpha}}(\mathbb{R})\models\phi_j(y,\eta),$$
then Player I has a winning strategy in the game obtained by playing the theory of $L_{\vartheta_{\alpha}}(\mathbb{R})$ in both subgames. This way, Player I is guaranteed to win the first subgame and Player II is forced to play the theory of $L_{\vartheta_{\alpha}}(\mathbb{R})$ in the second subgame, which will cost her the game. Similarly, if $\eta$ above exists and there is not a unique ordinal $\gamma<\vartheta_{\alpha}$ such that 
$$L_{\vartheta_{\alpha}}(\mathbb{R})\models\phi_i(x,\gamma),$$
then Player II has a winning strategy in the game.
Otherwise, let $\gamma$ and $\eta$ be as above. We claim that
\begin{equation}\label{eqLemma2TheoremBorel}
\gamma\leq\eta\text{ if, and only if, Player I has a winning strategy in the game.}
\end{equation}
If so, then the game defines a prewellordering of $\mathbb{R}$ of length greater than $\vartheta_{\alpha}^*$, in which the pairs of highest rank are those whose first coordinate fails the syntactic condition \eqref{LemmaBorelPWORule1}--\eqref{LemmaBorelPWORule2} and the pairs of second-highest rank are those
not defining an ordinal in $L_{\vartheta_\alpha}(\mathbb{R})$, i.e., those failing conditions \eqref{LemmaBorelPWORule4}--\eqref{LemmaBorelPWORule5}. The facts that condition \eqref{LemmaBorelPWORule1} is checked before condition \eqref{LemmaBorelPWORule2} and that condition \eqref{LemmaBorelPWORule4} is checked before \eqref{LemmaBorelPWORule5} will ensure that the game defines a reflexive relation.
Let us now proceed to the proof of equation \eqref{eqLemma2TheoremBorel}. 

Suppose $\gamma\leq\eta$. A winning  strategy for Player I is defined by setting $x_0 = x$ as appropriate and playing the theory of $L_{\vartheta_\alpha}(\mathbb{R})$ in the first subgame together with all reals played in the second subgame. In the second subgame, Player I also plays the theory of $L_{\vartheta_\alpha}(\mathbb{R})$. This strategy ensures that Player  I will win the first subgame, as in Lemma \ref{LemmaBorelVartheta}. Moreover, if Player II wishes to have any chance of winning the second subgame, she must also play the theory of $L_{\vartheta_\alpha}(\mathbb{R})$. If so, then this will define two models $N_I$ and $N_{II}$ as in the rules of the game, both wellfounded. Either there will be a real played in the first game and not in the second---in which case Player I wins---or both plays will contain the same reals. If so, the models $N_I$ and $N_{II}$ both embed into $L_{\vartheta_\alpha}(\mathbb{R})$ and are hence wellfounded, so their transitive collapses are of the forms $L_{\alpha_I}(\sigma)$ and $L_{\alpha_{II}}(\sigma)$, for some countable set of reals $\sigma$. Since each model satisfies the formulae:
\begin{enumerate}
\item $\kp$ + Separation,
\item \quotesleft $V_\alpha$ exists,'' and
\item 
\quotesleft for all $\eta \in \Ord$, $L_\eta(\mathbb{R})\not\models$ $\kp$ + Separation + \quoteleft $V_\alpha$ exists,' ''
\end{enumerate}
they must be equal. If Player II has not lost at this point, then we must have $y_0 = y$. The fact that $L_{\alpha_I}(\sigma)$ embeds elementarily into $L_{\vartheta_\alpha}(\mathbb{R})$ implies that the ordinal defined by $\phi_i$ and $x$ must be at most the ordinal defined by $\phi_j$ and $y$ in those models, as desired.

Conversely, suppose that Player I has a winning strategy in the game. We consider a run of the game by that strategy in which Player II plays the theory of $L_{\vartheta_\alpha}(\mathbb{R})$ in the first subgame, as well as in the second, where she also plays all reals from the first subgame. This ensures winning the second subgame and losing the first one. As before, this forces Player I to play the theory of $L_{\vartheta_\alpha}(\mathbb{R})$ in the first game, so the play defines two wellfounded models with the same reals and the same theory, whose transitive collapses must therefore be equal. Since the play must be won by Player I, the ordinal defined by $\phi_i$ and $x$ must be at most the ordinal defined by $\phi_j$ and $y$.

This proves the claim and, together with the remark following equation \eqref{eqLemma2TheoremBorel} shows that the relation $\preceq$ given by
\begin{equation*}
(i,x) \preceq (j,y) \text{ if, and only if, Player I has a winning strategy in $G(\alpha,i,x,j,y)$}
\end{equation*}
is a prewellordering of length greater than $\vartheta_\alpha^*$. A standard Skolem-hull argument shows that if $\beta<\alpha$, then $\vartheta_\beta<\vartheta_\alpha^*$.
Using the fact that the game from Lemma \ref{LemmaBorelVartheta} is Borel, one sees that for every $\alpha$, $G(\alpha,i,x,j,y)$ is uniformly in $\Delta^1_1(\alpha, \langle i,x \rangle, \langle j, y\rangle)$, which completes the proof.
\endproof

Before proving the theorem, we need two more lemmata.

\begin{lemma}\label{LemmaScaleVartheta}
Let $\omega+1< \alpha<\omega_1$. The pointclass $(\SIGMA^2_1)^{L_{\vartheta_\alpha}(\mathbb{R})}$ has the scale property.
\end{lemma}
\proof
Let $\gamma = \vartheta_\alpha$. Since $L_\gamma(\mathbb{R})$ is a model of $\Sigma_1$-separation, $\gamma$ is $\mathbb{R}$-nonprojectible, in the sense that there is no total function which is $(\Sigma_1)^{L_\gamma(\mathbb{R})}$ with parameters in $L_\gamma(\mathbb{R})$ and which maps $\gamma$ injectively into some $\beta<\gamma$. The usual argument for $L$ (see Barwise \cite[Theorem 6.8]{Ba75}) shows that $\gamma$ is recursively $\mathbb{R}$-inaccessible, i.e., $\mathbb{R}$-admissible and a limit of $\mathbb{R}$-admissibles. 

To prove the lemma, we first verify that
the proof of Martin-Steel \cite[Lemma 1]{MaSt08} goes through in $L_\gamma(\mathbb{R})$ and shows that
every $(\Game^\mathbb{R}\PI^1_1)^{L_\gamma(\mathbb{R})}$ set\footnote{In the choiceless context, $\Game^\mathbb{R}$ is to be understood in terms of quasi-strategies.} admits a closed game representation in $L_\gamma(\mathbb{R})$. To see this, first observe that every wellfounded tree $T$ in $L_\gamma(\mathbb{R})$ has a rank, since $\gamma$ is a limit of $\mathbb{R}$-admissibles. This observation ensures that the proof of the claim in the proof of \cite[Lemma 1]{MaSt08} holds in $L_\gamma(\mathbb{R})$ so that (using the notation from \cite{MaSt08}) Player I has a winning strategy in $G_x$ if, and only if, she has one in $G^*_{x,\xi}$ for some $\xi$. Now, each $G^*_{x,\xi}$ is a closed game and $G^*_{x,\xi} \in L_{\xi+1}(\mathbb{R})$. By standard arguments, if Player I has a winning strategy in $G^*_{x,\xi}$, she has one in $L_{\xi^*+1}(\mathbb{R})$ whenever $\xi<\xi^*$ and $\xi^*$ is $\mathbb{R}$-admissible. 
Clearly, $\Theta^{L_\gamma(\mathbb{R})}$ is also recursively $\mathbb{R}$-inaccessible and, by the L\"owenheim-Skolem theorem,
$$L_{\Theta^{L_\gamma(\mathbb{R})}}(\mathbb{R})\prec_1 L_\gamma(\mathbb{R}),$$
so it follows that
if $x\in\mathbb{R}$ and $\xi$ is least such that Player I has a winning strategy in $G^*_{x,\xi}$, then 
$$\xi<\Theta^{L_\gamma(\mathbb{R})}.$$
Denote this least $\xi$ by $\xi_x$ and let $\lambda$ be the supremum of all $\xi_x$ such that $\xi_x$ is defined; thus, $\lambda\leq \Theta^{L_\gamma(\mathbb{R})}$. This shows that the map $x\mapsto G^*_{x,\lambda}$ belongs to $L_\gamma(\mathbb{R})$. One finishes the argument as in \cite[Lemma 1]{MaSt08}.

The statement of \cite[Lemma 3]{MaSt08} also holds true in $L_\gamma(\mathbb{R})$, i.e.,
$$(\Game^\mathbb{R}\Pi^1_1)^{L_\gamma(\mathbb{R})} = \Sigma_1(L_\gamma(\mathbb{R}), \{\mathbb{R}\}).$$
Here, the proof needs no modifications other than weakening condition (2) to only require that the model satisfy (say) $\kp$, as opposed to $\zf^-$. This shows, using the proof of \cite[Theorem 4]{MaSt08}, that the pointclass 
$$\Sigma_1(L_\gamma(\mathbb{R}), \{\mathbb{R}\})$$
has the scale property. Now the proof of Solovay's basis theorem (e.g., the one in Section 2 of Koellner-Woodin \cite{KW10}, replacing their $T_0$ with some strong-enough theory that holds in $L_{\Theta^{L_\gamma(\mathbb{R})}}(\mathbb{R})$) shows that
$$\Sigma_1(L_\gamma(\mathbb{R}), \{\mathbb{R}\}) = (\Sigma^2_1)^{L_\gamma(\mathbb{R})},$$
from which the lemma follows.
\endproof

\begin{lemma}\label{LemmaShorteningBorelQuantifier}
Suppose that sets in $\Game^\mathbb{R}\DELTA^1_1$ are determined. Then, Borel games of length $\omega^2$ are determined.
\end{lemma}
\proof
This argument is a local version of the one in Blass \cite{Bl75}.
Consider the game on $\mathbb{R}$ of length $\omega$ in which two players alternate turns playing strategies for games of length $\omega$:
\begin{align*}
\begin{array}{c|ccccc}
I & \sigma_0 & & \sigma_1 & &\hdots\\
II & & \tau_0 & & \tau_1 & \hdots
\end{array}
\end{align*}
Player I wins if, and only if,
$$(\sigma_0*\tau_0, \sigma_1*\tau_1,\hdots) \in A,$$
where $\sigma*\tau$ denotes the result of playing the strategies $\sigma$ and $\tau$ against each other.

Call this game $G$. It is a Borel game (on $\mathbb{R}$), so it is determined. If it is Player I who wins this game, then Player I also wins the game of length $\omega^2$ with payoff $A$. Suppose that it is Player II who wins this game; we shall construct a winning strategy for the game of length $\omega^2$ with payoff $A$ for Player II.

Let $\Sigma$ denote the winning strategy for Player II in $G$. $\Sigma$ is essentially a set of reals. By the remark at the beginning of the proof of Lemma \ref{LemmaBorelVartheta}, we may assume 
$$\Sigma \in L_{\vartheta}(\mathbb{R}).$$
Recall that $\vartheta = \sup\{\vartheta_\alpha:\alpha<\omega_1\}$. By the minimality of $\vartheta_\alpha$ and the $\mathbb{R}$-stability of $\delta^2_1$, we have
$$\vartheta_\alpha < (\delta^2_1)^{\vartheta_{\alpha+1}}$$
for each $\alpha<\omega_1$. This implies by Lemma \ref{LemmaScaleVartheta} that $\mathcal{P}(\mathbb{R})\cap L_\vartheta(\mathbb{R})$ has the scale property. Since it is closed under quantification by $\mathbb{R}$, Theorem D of \cite{Mo71} implies that it has the uniformization property.

Denote by $A^*$ the set of all sequences $p$ of natural numbers of length $\omega\cdot n$, for some $n \in\mathbb{N}$, such that there is a sequence 
$$\vec \tau \in\mathbb{R}^n$$
such that $p$ results from applying $\Sigma$ to $\vec \tau$. Put
$$A^*_n = A^*\cap \mathbb{R}^n.$$
For any given $p \in A^*_1$, there might be many strategies $\tau$ for Player I witnessing $p \in A^*$; however, $A^*_1$ is definable from $\Sigma$, so it belongs to $L_\vartheta(\mathbb{R})$. Applying the uniformization property, one finds a function 
$$f_1:A^*_1\to\mathbb{R}$$
such that for each $p \in A^*_1$, $f_1(p)$ is a strategy $\tau$ witnessing $p \in A^*$ and, moreover, $f_1 \in L_\vartheta(\mathbb{R})$. Inductively, suppose $f_n$ has been defined. For any given $p \in A^*_{n+1}$, there might be many strategies $\tau$ such that $f_n(p\upharpoonright\omega\cdot n)^\frown \tau$ witnesses $p \in A^*$, or none at all; however, $A^*_{n+1}$ is definable from $\Sigma$, so it belongs to $L_\vartheta(\mathbb{R})$. Applying the uniformization property, one finds a function 
$$f_{n+1}:A^*_{n+1}\to\mathbb{R}^{n+1}$$
with $f_{n+1} \in L_\vartheta(\mathbb{R})$ and such that the following hold for each $p \in A^*_{n+1}$ for which such a $\tau$ exists:
\begin{enumerate}
\item $f_{n+1}(p)$ is a sequence $\vec \tau$ witnessing $p \in A^*$; and
\item $f_{n+1}(p)\upharpoonright m = f_m(p\upharpoonright \omega\cdot m)$
for every $m\leq n+1$.
\end{enumerate}
Let us also define a partial function $f_\omega: \mathbb{R}^{\mathbb{N}} \to \mathbb{R}^\mathbb{N}$ by
$$f_\omega(p) = \bigcup_{n\in\mathbb{N}} f_{n}(p\upharpoonright  n),$$
if the right-hand side is defined, and put
$$f = \bigcup_{n\leq\omega}f_n.$$

Let us call a sequence $p$ of natural numbers of length $\omega\cdot n$, with $n\in\mathbb{N}$, \emph{promising} if  for every $m \leq n$, $p\upharpoonright \omega\cdot m \in A^*$ and this is witnessed by $f(p)\upharpoonright m$.
Given a promising sequence $p$, we shall construct a strategy $\sigma^*$ for Player II (for a game of length $\omega$ on $\mathbb{N}$) such that for every
$x\in\mathbb{R}$, $p^\frown(\sigma^*(x))$ is promising. 
Thus, if $p^*$ is a run of a game of length $\omega^2$ all of whose initial segments of limit length were obtained in this way, then it is a play all of whose initial segments are promising. This means that
$$p^* = f(p^*) * \Sigma\big(f(p^*)\big),$$
which implies that $p^* \in A$.

Now, let $p \in\mathbb{R}^n$ be promising and consider the game $G^*$ on $\mathbb{N}$ of length $\omega$:
\begin{align*}
\begin{array}{c|ccccc}
I & x(0) & & x(2) & &\hdots\\
II & & x(1) & & x(3) & \hdots
\end{array}
\end{align*}
Player II wins if, and only if, $p^\frown x$ is promising. We first show that Player I cannot have a winning strategy. Suppose towards a contradiction that $\tau$ is a winning strategy for Player I.  Thus, if Player I begins a run of $G$ by playing
$$f(p)^\frown\tau,$$
during the first $n+1$ moves and Player II plays by $\Sigma$, then
she will begin by responding with $n$ moves according to $\Sigma$, resulting in a sequence of strategies $\sigma_0,\hdots, \sigma_{n-1}$ so that
$$p = \big(f(p)(0)*\sigma_0, \hdots, f(p)(n-1)*\sigma_{n-1}\big)$$
and continue responding with one last move by $\Sigma$, say $\sigma$. Put $x = \tau*\sigma$. Since the play of $G$ just described is consistent with $\Sigma$, we have $p^\frown x \in A^*$. Although it need not be the case that $f(p^\frown x)(n) = \tau$,
the construction of $f$ ensures that $f(p^\frown x)$ is a sequence $\vec\tau$ witnessing $p^\frown x\in A^*$ and $f(p^\frown x)\upharpoonright n = f(p)$. Thus, the run of $G^*$ in which Player I plays by $\tau$ and Player II plays by $\sigma$ results in a play $x$ such that $p^\frown x$ is promising, which is a contradiction. Hence, Player I cannot have a winning strategy in $G^*$. We will show that the payoff set of $G^*$ belongs to $\Game^\mathbb{R}\DELTA^1_1$; this will imply that it is determined and thus that Player II has a winning strategy, which will finish the proof.

The payoff set of $G^*$ (for Player II) consists of all $x\in\mathbb{R}$ such that $p^\frown x$ is promising. A sequence being promising is defined in terms of $\Sigma$ and $f$. Now, each $f_n$ belongs to $L_\vartheta(\mathbb{R})$ and, since
$$\cof(\vartheta) = \omega_1,$$
there is some $\alpha<\omega_1$ such that  each $f_n$ belongs to $L_{\vartheta_\alpha}(\mathbb{R})$. Since 
$$\big(L_{\Theta^{L_{\vartheta_\alpha}(\mathbb{R})}}(\mathbb{R})\big)^\omega\subset L_{\vartheta_\alpha}(\mathbb{R}),$$
we have $f \in L_{\vartheta_\alpha}(\mathbb{R})$. By choosing $\alpha$ large enough so that $\Sigma \in L_{\vartheta_\alpha}(\mathbb{R})$, we ensure that
$$\big\{p \in \mathbb{R}^n: p \text{ is promising and } n\in\mathbb{N}\big\} \in L_{\vartheta_\alpha}(\mathbb{R}).$$
By Lemma \ref{LemmaBorelPWO}, $\vartheta\leq\gamma^1_1$, so that, by Lemma, \ref{TheoremDG}, the set in the displayed equation belongs to $\Game^\mathbb{R}\DELTA^1_1$, which completes the proof.
\endproof

Woodin has shown that $\mathbb{R}^\sharp$ exists and $\ad$ holds in $L(\mathbb{R})$ if, and only if, there is a countably iterable proper-class model of $\zfc$ with infinitely many Woodin cardinals. Neeman \cite{Ne04} has shown that if such a model exists, then all analytic games of length $\omega^2$ are determined (incidentally, the converse follows from a theorem of Trang \cite{Tr13} and a theorem of Martin-Steel \cite{MaSt08}, after appealing to Woodin's theorem). 
Let us remark that the argument just given yields a direct proof of analytic determinacy for games of length $\omega^2$ which involves no inner model theory.

\begin{theorem}(Neeman-Woodin)
Suppose that $\mathbb{R}^\sharp$ exists and $L(\mathbb{R})\models\ad$. Then, analytic games of length $\omega^2$ are determined.
\end{theorem}
\proof
Given an analytic set $A$, let $G$ be the game defined as in the proof of Lemma \ref{LemmaShorteningBorelQuantifier}. It is an analytic game on $\mathbb{R}$ and, if won by Player I, then the game of length $\omega^2$ with payoff $A$ is also won by Player I. 
Since $\mathbb{R}^\sharp$ exists, then analytic games on $\mathbb{R}$ are determined, by the argument in Martin \cite{Ma70}, so if Player I does not win $A$, then Player II does. 

Since $\ad$ holds in $L(\mathbb{R})$, all sets in $\Game^\mathbb{R}\PI^1_1$ are determined. Moreover, the pointclass $\Game^\mathbb{R}\PI^1_1$ is correctly computed in $L(\mathbb{R})$, so if Player II wins an analytic game on $\mathbb{R}$ (this is a $\Game^\mathbb{R}\PI^1_1$ fact), then she has a winning strategy which is $(\DELTA^2_1)^{L(\mathbb{R})}$. Hence, the game $G^*$ defined as in the proof of the preceding lemma belongs to $L(\mathbb{R})$ and is thus determined. The rest of the argument is as before.
\endproof

Theorem \ref{mainborel} now follows: if Borel games of length $\omega^2$ are determined, then all sets in 
$\Game^\mathbb{R}\DELTA^1_1$
are determined (in the usual sense). By Lemma \ref{LemmaBorelPWO}, $\vartheta\leq\gamma^1_1$, so 
\begin{equation}\label{eqBorel1}
L_\vartheta(\mathbb{R})\models\ad
\end{equation}
by Lemma \ref{TheoremDG}. Conversely, suppose \eqref{eqBorel1} holds. By Lemma \ref{LemmaBorelVartheta} and the fact that the existence of winning strategies for games on $\mathbb{R}$ is upwards absolute from initial segments of $L(\mathbb{R})$ (since they contain all possible plays of the game), we have 
\begin{equation}\label{eqBorel2}
\Game^\mathbb{R}\DELTA^1_1 \subset L_\vartheta(\mathbb{R}).
\end{equation}
Thus, sets in $\Game^\mathbb{R}\DELTA^1_1$ are determined (in the usual sense). By the previous lemma, this implies that games of length $\omega^2$ with Borel payoff are determined.
Theorem \ref{mainborel} is thus proved.

\begin{corollary}
$\gamma^1_1 = \vartheta$.
\end{corollary}
\proof
Immediate from Lemma \ref{TheoremDG}, Lemma \ref{LemmaBorelPWO}, and the observation that Lemma \ref{LemmaBorelVartheta} implies equation \eqref{eqBorel2}. Note that the three lemmata were proved under no determinacy hypotheses.
\endproof

\section{Proof of Lemma \ref{TheoremDG}}
\label{SectDG}
Let $\alpha<\gamma^1_1$ be an ordinal; we show
$$\mathcal{P}(\mathbb{R})\cap L_{\alpha+1}(\mathbb{R})\subset\Game^\mathbb{R}\DELTA^1_1.$$
The proof consists in defining a game that determines membership in a set definable over $L_\alpha(\mathbb{R})$; we call it the \emph{definability game}. Readers familiar with Tanaka's \cite{Ta90} proof of Steel's \cite{St77} classical result on the equivalence between arithmetical transfinite recursion and clopen determinacy should find similarities.

Fix some prewellordering $\preceq$ of $\mathbb{R}$ in $\Game^\mathbb{R}\DELTA^1_1$ of length $\alpha+1$ or greater. To ease notation, assume that $\preceq$ is defined without real parameters.
The definability game will involve playing \emph{representations} of elements of $L_{\alpha}(\mathbb{R})$. These are defined inductively: pairs $(0, x)$, where $x\in\mathbb{R}$, are representations for $x$. Let $\beta+1\leq\alpha$, $b \in L_\beta(\mathbb{R})$ and
$$B = \{x\in L_{\beta}(\mathbb{R}): L_{\beta}(\mathbb{R})\models\hat\phi(x,b)\}$$
be an element of $L_{\beta+1}(\mathbb{R})$. If $x_b$ codes a representation $(\psi, y, y_b)$ of $b$ and $|y|_\preceq < \beta$, then the triple
\begin{equation}\label{eqRepresentation}
(\hat\phi,x_\beta, x_b),
\end{equation}
is a representation of $B$, whenever $x_\beta$ is any real of rank $\beta$ in $\preceq$.\footnote{Strictly speaking, we should allow $x_b$ to code finite tuples of representations, but we shall generally abuse notation by assuming that these tuples have length $1$.}

\begin{lemma}\label{LemmaRepresentationsExist}
Let $\alpha$ and $\preceq$ be as above. Then, every set in $L_\alpha(\mathbb{R})$ has a representation.
\end{lemma}
\proof
By induction on ordinals $\beta\leq\alpha$.
\endproof

Let us also inductively define a \emph{proto-representation} to be either a pair $(0,x)$ as above, or a triple as in \eqref{eqRepresentation}
whose first coordinate is a formula and third coordinate is a proto-representation, i.e., one drops the requirement on the $\preceq$-rank of the second coordinate.

Being a proto-representation is easy to check---the class of all proto-representations is $\Sigma^0_1$. Being a genuine representation, on the other hand, is more complicated, as, if $(\hat \phi, x_\beta, x_b)$ is as in \eqref{eqRepresentation}, then one needs to check that $x_b$ in fact codes a representation for an element of $L_\beta(\mathbb{R})$.

In principle, one can have highly inefficient representations in which one defines simple objects in a complicated way. We will denote by $\rep_\alpha$ the set of codes of all representations of the form $(0,x)$ and of the form $(\hat\phi, x_\beta, x_b)$, where $x_\beta$ has rank ${<}\alpha$ in $\preceq$. Belonging to $\rep_\alpha$ does not guarantee that a representation is optimal, but it guarantees that, in a sense, it is not \quotesleft more complicated than $\alpha$.''

\begin{lemma}\label{LemmaRepCond}
Let $\alpha<\gamma^1_1$, $\preceq$, and $\rep_\alpha$ be as above. If $x_\alpha$ has rank $\alpha$ in $\preceq$, then $\rep_\alpha \in \Game^\mathbb{R}\Delta^1_1(x_\alpha)$.
\end{lemma}
\proof
It is easy to determine whether a real number codes a representation of the form $(0,y)$, so let us consider only proto-representations which are triples.

For a triple $(\psi, x, y)$ to be a representation of an element of $L_\alpha(\mathbb{R})$ in $\rep_\alpha$, it is necessary and sufficient that the following two conditions hold:
\begin{enumerate}
\item \label{LemmaRepCond1} $x \prec x_\alpha$,
\item \label{LemmaRepCond2} $y$ codes a representation $(\hat\psi, x_\gamma, x_c)$ and $x_\gamma\prec x$.
\end{enumerate}
Given $a \in \mathbb{R}$, we define a game $G(a)$ on reals with payoff in $\Delta^1_1(x_\alpha)$. In this game, Player I claims that $a$ codes a representation, and Player II claims otherwise. 

If $a$ does not code a triple, the game ends immediately and Player I loses. Otherwise, suppose $a$ codes a triple $(\psi, x, y)$ as above.
Player II must begin by challenging one of the two conditions \eqref{LemmaRepCond1} or \eqref{LemmaRepCond2} above. If Player II decides to challenge condition \eqref{LemmaRepCond1}, then they must play the game with payoff in $\Delta^1_1(x,x_\alpha)$ determining whether 
$$x \prec x_\alpha;$$
if so, the winner of this subgame is the winner of $G$. If Player II decides to challenge condition \eqref{LemmaRepCond2}, then Player II must provide one of the following three reasons:
\begin{enumerate}
\item $y$ does not code a triple nor a pair $(0,\hat y)$, where $\hat y \in \mathbb{R}$;
\item $y$ codes a triple $(\hat\psi, x_\gamma, x_c)$, but $x\preceq x_\gamma$;
\item $y$ codes a triple $(\hat\psi, x_\gamma, x_c)$, but $x_c$ is not a representation.
\end{enumerate}
If Player II claims that $y$ does not code a triple or a pair $(0,\hat y)$, then the game ends and Player I wins if, and only if, $y$ codes a triple or a pair $(0,\hat y)$. If Player II claims that $x\preceq x_\gamma$, then, as above, they must play the game with payoff in $\Delta(x,x_\gamma)$ determining whether
$$x_\gamma \prec x.$$
Finally, if Player II claims that $y$ codes a triple $(\hat\psi, x_\gamma, x_c)$, but $x_c$ is not a representation, then she again must provide a reason, and so on. 
If Player II correctly identifies the problems with Player I's purported representations throughout the game, then she must initiate a subgame after finitely many turns (because $\prec$ is wellfounded). Hence, we add as a rule that if the game finishes after infinitely many turns without Player II having initiated a subgame, she loses. We note that, although different subgames can be initiated by Player II, they are all Borel, and in fact they are so in a uniform way, i.e., they are all instances of the same game (namely, the one for determining whether $u \prec v$), and only their parameters vary. This uniformity is necessary for ensuring that the payoff of $G(a)$ is also Borel.

Now, Player I wins a run $p$ of the game if, and only if, one of the following holds:
\begin{enumerate}
\item One of the players violates a rule and the first one to do so is Player II; 
\item The players obey the rules and for all $n$, Player II does not initiate a subgame in the $n$th turn; 
\item The players obey the rules, there is a least $n$ such that Player II initiates a subgame in the $n$th turn, and $p\upharpoonright(n,\infty)$ is a win for Player I in this game.
\end{enumerate}
The first condition is $\Sigma^0_1(a,p, x_\alpha)$, the second one is $\Pi^0_1(a,p)$, and the third one is $\Delta^1_1(a,p)$. Thus, the winning condition is in $\Delta^1_1(a,p,x_\alpha)$. Moreover, if $(\psi, x, y)$ is a representation of an element of $L_\alpha(\mathbb{R})$, then any challenge raised by Player II can easily be overcome by Player I, so she must have a winning strategy. Similarly, if $(\psi, x, y)$ is not a representation of an element of $L_\alpha(\mathbb{R})$, then Player II can identify the reason why this is the case and pose a challenge Player I cannot overcome. Therefore, Player I has a winning strategy in this game if, and only if, $a$ indeed codes a representation of an element of $L_\alpha(\mathbb{R})$ in $\rep_\alpha$. We have shown that the set of codes of representations of elements of $L_\alpha(\mathbb{R})$ is the set of all $a$ for which Player I has a winning strategy in the Borel game $G(a)$, as was to be shown.
\endproof

\begin{remark}
The fact that all auxilary games were  Borel in a uniform way was crucial in the proof of Lemma \ref{LemmaRepCond}. In fact, the set $\bigcup_{\alpha<\gamma^1_1}Rep_\alpha$ is not in $\Game^\mathbb{R}\DELTA^1_1$. \qed\\
\end{remark}

\begin{definition}
The \emph{order} of a proto-representation is defined inductively: the order of a proto-representation of the form $(0,x)$ is $0$; the order of a proto-representation of the form $(\phi,x_\beta, x_b)$ is the maximum of $|x_\beta|_\prec$ and the order of $x_b$.
\end{definition}
Thus, if $x \in \rep_\alpha$, then the order of $x$ is strictly less than $\alpha$.\\

We now proceed to the definability game and, with it, to the proof of Lemma \ref{TheoremDG}. The idea of the game is very simple: two players argue whether a formula holds of some sets in a given level of the $L(\mathbb{R})$-hierarchy. However, its description is lengthy. This is due to two complications that arise: the first one is that we need the payoff of the game to be $\DELTA^1_1$. This requires us speaking about representations of sets instead of sets directly. The second one is that many rules need to be put into place to ensure that the players are honest in their moves.

Fix some $n\in\mathbb{N}$, a parameter $x_A\in L_\alpha(\mathbb{R})$, and a set 
$$A \in \Sigma_n(L_{\alpha}(\mathbb{R}),x_A),$$
say, $A = \{y\in L_{\alpha}(\mathbb{R}): L_{\alpha}(\mathbb{R})\models\phi(y,x_A)\}.$
In the definability game, Player I attempts to show that a given $y \in L_{\gamma^1_1}(\mathbb{R})$ belongs to the set $A$ above, i.e., that
\begin{equation}\label{eqDefGame1}
L_{\alpha}(\mathbb{R})\models\phi(y,x_A)
\end{equation}
We assume that the scopes of all negations in $\phi$ (and in all formulae appearing below) are atomic formulae. We also assume that no implications appear in $\phi$. The game will involve Players I and II playing purported representations of sets in $L_\alpha(\mathbb{R})$. 
Fix some representation $(\phi,x_\alpha, x_a)$ for $A$ and some real $x^*$ coding it. Let $a\in\mathbb{R}$; we define the game $D_A(a)$. If $a$ does not code a pair $(0,y)$ or a triple $(\phi_y, x_y, x_c)$, the game ends immediately and Player I loses. Otherwise, the game begins. We will describe the rules in the case that $a$ codes a triple; the rules in the case that $a$ codes a pair $(0,y)$ are analogous; soon, we will restrict our attention to reals $a$ of order ${<}\alpha$ (as the reader will soon realize, this restriction postpones a certain complication in the definition of the game). 
Player I's goal is to convince Player II that the triple coded by $a$ is a representation for some $y$ such that equation \eqref{eqDefGame1} is satisfied. The game is defined from $x^*$ and $a$ (as well as the parameters defining $\preceq$, which we assumed for simplicity to be non-existent). 

The rules of the game give the players opportunities to initiate subgames in which they challenge the other's moves. If this happens at any point, the definability game will end and the winner will be declared to be the winner of the subgame. We do not need to speak about winning strategies for the subgame---we can simply ask the players to play it. All subgames will have Borel payoff and real moves; this allows the players to check during the game whether a given $\Game^\mathbb{R}\DELTA^1_1$ fact holds.\\

At the beginning of of $D_A(a)$, Player II has the oportunity to claim that
\begin{align*}
\text{$(\phi_y, x_y, x_c)$ is not a representation of an element of $L_\alpha(\mathbb{R})$,}
\end{align*}
i.e., that $(\phi_y, x_y, x_c)\not\in\rep_\alpha$.
If this happens, then the game is over and Player I wins if, and only if, the challenged object is indeed a representation; otherwise Player II wins. If Player II does not claim either of the two statements above, the game proceeds.

At the beginning of turn $k$, Players I and II have defined a proto-representation $(\phi_k, x_{\beta^k}, x_{b^k})$ for some set $B^k$ and finitely many proto-representations $(\psi_{w^k_i}, x_{w^k_i}, x_{d^k_i})$ for some sets $w^k_i$. Player I claims that $\vec w^k \in B^k$, i.e., that
$$L_{\beta^k}(\mathbb{R})\models\phi_k(\vec w^k,b^k);$$
Player II claims otherwise.
Below, we write $\vec w = \vec w^k$ and $b = b^k$ for notational simplicity; as before, we might abuse notation by assuming $\vec w$ is a tuple of length $1$ and simply writing $w$.
We define the \emph{order} of the game at turn $k$ to be the maximum of the orders of $w$ and $b$.
Let us describe the rules of the game by cases:

\begin{enumerate}
\item If $\phi_k$ is of the form $\psi\vee \chi$, then Player I selects either $\psi$ or $\chi$. We set $\phi_{k+1}$ equal to that choice and leave the rest of the objects unchanged.
\item If $\phi_k$ is of the form $\psi\wedge \chi$, then Player II selects either $\psi$ or $\chi$. We set $\phi_{k+1}$ equal to that choice and leave the rest of the objects unchanged.
\item If $\phi_k(w, b)$ is of the form $\exists z\, \psi(w, b, z)$, then Player I must play a proto-representation $x_z$ which is allegedly a representation in $\rep_{\beta^k}$ of some $z\in L_{\beta^k}(\mathbb{R})$ such that
$$L_{\beta^k}(\mathbb{R})\models \psi(w,b,z).$$
At this point, Player II has the opportunity to object to the fact that $x_z$ belongs to $\rep_{\beta^k}$. If so, the game ends and Player I wins if, and only if, $x_z \in \rep_{\beta^k}$. If Player II does not challenge Player I's move, then the game continues with $\phi_{k+1} = \psi$ and $z$ added to $\vec w$; the rest of the objects remain unchanged.
\item If $\phi_k(w, b)$ is of the form $\forall z\, \psi(w, b, z)$, then Player II must play a proto-representation $x_z$ which is allegedly a representation in $\rep_{\beta^k}$ of of some $z\in L_{\beta^k}(\mathbb{R})$ such that
$$L_{\beta^k}(\mathbb{R})\not\models \psi(w,b,z).$$
As before, Player I has the opportunity to object to the fact that $x_z$ belongs to $\rep_{\beta^k}$. If so, the game ends and Player I wins if, and only if, $x_z\not\in\rep_{\beta^k}$. If Player I does not challenge Player II's move, then the game continues with $\phi_{k+1} = \psi$ and $z$ added to $\vec w$; the rest of the objects remain unchanged.
\item \label{condDefGame5} If $\phi_k$ is a non-negated atomic formula, then it is of the form 
$u \in v$, where $u$ and $v$ are any of $b$ or some $w_i$; both $u$ and $v$ have associated proto-representations that have been played in one of the previous turns, or given by the initial data. There are three subcases that we need to distinguish: 
\begin{enumerate}
\item Suppose that the proto-representation associated to $v$ is of the form $(0, x)$, for some $x\in\mathbb{R}$, and the proto-representation associated to $u$ is of the form $(0, y)$, for some $y\in\mathbb{R}$; then the game ends. Player I wins if, and only if, $y \in x$; otherwise, Player II wins. 
\item \label{condDefGame5b} Suppose that the proto-representation associated to $v$ is of the form $(0, x)$, for some $x\in\mathbb{R}$, but the one associated to $u$ is of the form
$$(\phi_u, x_{\eta_u}, x_{u}).$$
The problem is that---in principle---$(\phi_u, x_{\eta_u}, x_{u})$ might be a complicated (true) representation of a simple set.\footnote{Originally introduced into the game by either player.} The game proceeds as follows: since the proto-representation associated to $v$ is of the form $(0, x)$, it is a true representation of some real number. 
Player I must thus play a natural number $n \in x$ and claim that $(\phi_u, x_{\eta_u}, x_u)$ is a representation of $n$, i.e., that, there is $u' \in L_{\eta_u}(\mathbb{R})$ represented by $x_u$ such that
\begin{equation*}
L_{\eta_u}(\mathbb{R})\models \phi_u(u',m) \text{ if, and only if, } m<n.
\end{equation*}
Player II must object by playing some $m_0 \in\mathbb{N}$ such that one of the following holds:
\begin{enumerate}
\item $m_0 < n$ but $L_{\eta_u}(\mathbb{R})\not\models \phi_u(u',m_0)$; or
\item $ n < m_0$ and $L_{\eta_u}(\mathbb{R})\models \phi_u(u',m_0)$.
\end{enumerate}

During the remainder of the game, the two players must determine whether 
$$L_{\eta_u}(\mathbb{R})\models \phi_u(u',m_0).$$
We know how to do this: namely, we set $B^{k+1}$ equal to
$$\big\{m\in\mathbb{N}: L_{\eta_u}(\mathbb{R})\models \phi_u(u',m)\big\}$$
or
$$\big\{m\in\mathbb{N}: L_{\eta_u}(\mathbb{R})\models \lnot\phi_u(u',m)\big\},$$
according to which objection was raised by Player II, and work with the canonical representation for $B^{k+1}$; and $\vec w^{k+1} = (0, m_0)$.
\item Suppose that $v$ has an associated proto-representation of the form
$$(\chi, x_\eta, x_{a_v}).$$
This case is similar to the previous one: Player I claims $u \in v$, i.e., 
$$L_{\eta}(\mathbb{R})\models \chi(u,a_v).$$
The problem (as before) is that the proto-representation
$$(\phi_u, x_{\eta_u}, x_{u})$$
associated to $u$
might be unnecessarily complicated. We ask Player I to decide whether
$x_{\eta_u}\prec x_{\eta}$.
To ensure that Player I tells the truth, we then ask Player II to decide whether 
$x_{\eta_u}\prec x_{\eta}$.
If the two players disagree on this, the game ends, with the winner being declared as the winner of the subgame that determines whether $x_{\eta_u} \prec x_{\eta}$.

Otherwise, suppose first that the players both believe that $x_{\eta_u}\prec x_{\eta}$. 
At this point, Player II can claim that the proto-representation $(\phi_u, x_{\eta_u}, x_{u})$ associated to $u$ does not belong to $Rep_\eta$. If so, then---as above---the game ends, and the winner of the game is the winner of the game given by Lemma \ref{LemmaRepCond} applied to $x_\eta$ and the proto-representation associated to $u$.\footnote{Note that this proto-representation could have been played by Player II at an earlier turn; this does not matter.} If Player II chooses not to make that claim, the game continues with $\phi_{k+1} = \chi$, $x_{\beta^{k+1}} = x_\eta$, $x_{b^{k+1}} = x_{a_v}$, and $\vec w^{k+1}=u$.

If the players both believe that $x_{\eta} \preceq x_{\eta_u}$, then we need to argue as above. If $u$ truly is an element of $v$, then $u$ must have a representation of order strictly less than $\eta$. Player I must play such a representation, say, $(\phi_{u'}, x_{\eta_{u'}}, x_{u'})$. As usual, Player II may object by claiming that $(\phi_{u'}, x_{\eta_{u'}}, x_{u'}) \not \in \rep_\eta$, in which case the game ends and the winner is declared according as whether $(\phi_{u'}, x_{\eta_{u'}}, x_{u'}) \in \rep_\eta$ or not; or by claiming that $(\phi_{u'}, x_{\eta_{u'}}, x_{u'})$ is not in fact a representation of the same object as $(\phi_u, x_{\eta_u}, x_{u})$. Let $u'$ be the object represented by $(\phi_{u'}, x_{\eta_{u'}}, x_{u'})$ (if any). This is now treated much like the case \eqref{condDefGame5b}. Player II must claim one of the following:
\begin{enumerate}
\item $u' \not\subset u$; or
\item $u\not\subset u'$.
\end{enumerate} 
Assume without loss of generality that it is the first alternative that is claimed. Player II must play a proto-representation $(\phi_{v'}, x_{\eta_{v'}}, x_{v'})$ for a set $v'$ witnessing that $u' \neq u$. Player I can make the usual objection that $(\phi_{v'}, x_{\eta_{v'}}, x_{v'})\not \in \rep_{\eta'_u}$, in which case the game proceeds as usual. Otherwise, Player I may claim that $v' \not \in u'$ or that $v' \in u$. In the former case, the game proceeds with the proto-representation $(\phi_{u'}, x_{\eta_{u'}}, x_{u'})$ for $u' = B^{k+1}$ and $w^{k+1}$ equal to the proto-representation $(\phi_{v'}, x_{\eta_{v'}}, x_{v'})$. The other case is similar.
\end{enumerate}

\item \label{condDefGame6} Finally, if $\phi_k$ is a negated atomic formula, then it is of the form 
$u \not\in v$, where $u$ and $v$ are any of $b$ or some $w_i$; both $u$ and $v$ have associated proto-representations that have been played in one of the previous turns. One distinguishes three subcases and proceeds as in case \eqref{condDefGame5}.
\end{enumerate}

The hope is that
both players play only real representations when they should play proto-representations and that the initial triples are also real representations. 
Let us say that a (partial) play in which these conditions are satisfied is \emph{honest}.

A key remark is that, if a play $p$ is honest, then the order of $p$ decreases after each time the players find themselves in one of situations \eqref{condDefGame5} or \eqref{condDefGame6} (assuming the game does not end).
This implies that an honest play in which both players
accept each other's moves will end after finitely many turns by one of the winning conditions in clauses \eqref{condDefGame5}--\eqref{condDefGame6} and the winner will be defined then.

If our expectations on the honesty of the players are not realized, the game will end in one of two possible ways: by means of a subgame initiated by a player challenging the other player's choice, in which case the winner of the definability game will be the winner of the subgame; or perhaps after infinitely many turns, if a player made improper moves and was not challenged by the opponent (or perhaps if one of the initial triples was not a real representation and this was not noticed by Player II). In this latter case, we declare Player I as the winner.\footnote{This choice will not be very consequential. Such an outcome will not occur in any case of interest, but it needs to be considered to make the game zero-sum.} 

There are three types of subgames that can be initiated: the one given by Lemma \ref{LemmaRepCond}, its dual---that in which the roles of Player I and Player II are reversed---, and the game associated to the prewellordering $\prec$. Games of each of those three types are uniformly in $\DELTA^1_1$, and defined using the initial parameters, as well as parameters played during the (main) game.
 
Thus, Player I wins a run $p$ of the game if, and only if, $a$ codes a triple and
\begin{enumerate}
\item There is some $k$ such that a player breaks a rule during turn $k$ and the first player to do so is Player II;
\item There is some $k$ such that a player initiates a challenge during turn $k$ and $p\upharpoonright(k,\infty)$ is a win for Player I in the subgame;
\item The game ends at a finite stage $k$ by Player I fulfilling the condition in case \eqref{condDefGame5} or case \eqref{condDefGame6} above; or
\item After infinitely many steps, no player has won or initiated a challenge.
\end{enumerate}
The first condition is $\Sigma^0_1(x^*,a,p)$. The second condition is $\Delta^1_1(x^*,p)$. Lastly, the third condition is $\Sigma^0_1(p)$ and the fourth condition is $\Pi^0_1(p)$. Hence, the payoff set is in $\Delta^1_1(x^*,a,p)$. 
\begin{lemma}
Player I has a winning strategy in $D_A(a)$ if, and only if, $a$ codes a representation in $\rep_\alpha$ of some set $y \in L_\alpha(\mathbb{R})$ such that
\begin{equation}\label{eqDefGameLemma}
L_\alpha(\mathbb{R})\models\phi(y,x_A).
\end{equation}
\end{lemma}
\proof
Clearly it is necessary that $a$ code a representation in $\rep_\alpha$ for some set $y \in L_\alpha(\mathbb{R})$ in order for Player I to have a winning strategy.
Suppose it does and that \eqref{eqDefGameLemma} holds.
The strategy for Player I is simple---essentially, she will always tell the truth. This will ensure that the intuitive interpretation of the game in terms of Player I claiming that a formula holds in an initial segment of $L(\mathbb{R})$ and Player II claiming otherwise is accurate.

If Player I's assertions are ever challenged by Player II, the subgame can be won as in Lemma \ref{LemmaRepCond}. Otherwise, at a given turn $k$, if $(\phi_k, x_{\beta^k},x_{b^k})$ and $(\psi_{w^k}, x_{w^k},x_{d^k})$ are as in the definition of the game, then they will be true representations. If $\phi_k$ is a disjunction, then one of the disjuncts will hold in $L_{\beta^k}(\mathbb{R})$ and Player I will choose it. If $\phi_k$ is a conjunction, then both conjuncts will hold in $L_{\beta^k}(\mathbb{R})$ and so it will not matter which one Player II chooses. If $\phi_k$ is of the form $\exists z\, \psi$, then there really must be a witness to $\psi$ in $L_{\beta^k}(\mathbb{R})$, and Player I will play a true representation of this witness, say, of minimal order. If $\phi_k$ is of the form $\forall z\, \psi$ and Player II plays a proto-representation which is not a representation of an element of $L_{\beta^k}(\mathbb{R})$, or has order higher than that permissible by the rules of the game, then Player I will challenge Player II's move and win the game as in Lemma \ref{LemmaRepCond}. Otherwise, Player II will play a representation of an element of $L_{\beta^k}(\mathbb{R})$ and $\psi$ will hold of this element. 
If $\phi_k$ is atomic or negated atomic, then it will be true, so that either the game will end and be won by Player I or continue according to the form of the representations relevant to the turn of the game the players are in. As remarked earlier, the order of every honest play of the game decreases after each time the players find themselves in one of situations \eqref{condDefGame5} or \eqref{condDefGame6}.
Hence, the strategy just described  ensures that either a player will initiate a challenge that will end by Player I winning, or that the game will end after finitely many turns and be won by Player I.

Similarly, if Player I has a winning strategy, then \eqref{eqDefGameLemma} must hold, for otherwise essentially the same strategy for Player II will be a winning strategy. The point is that if Player II plays properly, every play will end either after finitely many turns or with a subgame.
\endproof

Thus, given a set of reals $A \in L_{\alpha+1}(\mathbb{R})$
and a code $x^*$ of a representation for $A$ (which exists by Lemma \ref{LemmaRepresentationsExist}), one has
$$A = \Big\{y\in\mathbb{R}: \text{ Player I has a winning strategy in } D_A(\langle 0, y\rangle)\Big\}.$$
Since $D_A(\langle 0, y\rangle)$ is in $\Delta^1_1(x^*, \langle 0, y\rangle)$, we have
$$\mathcal{P}(\mathbb{R})\cap L_{\alpha+1}(\mathbb{R})\subset\Game^\mathbb{R}\DELTA^1_1.$$
Therefore, the proof of Lemma \ref{TheoremDG} is complete.

\section{Derived model theorems} \label{BDSectDMT}

In this section, we prove two lemmata that tie determinacy in $L_\vartheta(\mathbb{R})$ to the existence of countably iterable extender models of Zermelo set theory with infinitely many Woodin cardinals. The first lemma is proved just like its analogue for $L(\mathbb{R})$. The second lemma is also similar to its analogue for $L(\mathbb{R})$ but has some differences.
We will assume some familiarity with the theory of extender models, e.g., as presented in Mitchell-Steel \cite{MS94} or Steel \cite{St10}.

We note that every interval of the form $[(\delta^2_1)^{L_{\vartheta_\alpha}(\mathbb{R})},\vartheta_\alpha]$ is a gap (cf. the proof of Lemma \ref{LemmaScaleVartheta} above).

\begin{lemma}\label{LemmaBorelDM}
Suppose that for every $\alpha<\omega_1$, there is a countable, countably iterable extender model of
$$\mathsf{Z} + \text{\quotesleft$V_{\lambda+\alpha}$ exists and $\lambda$ is a limit of Woodin cardinals.''}$$
Then $L_{\vartheta}(\mathbb{R})\models\ad$.
\end{lemma}
\proof
Let $\alpha<\omega_1$ be large. We show that
\begin{equation}\label{LemmaBorelDMMaineq}
L_{\vartheta_\alpha}(\mathbb{R})\models\ad.
\end{equation}

Let $M$ be an extender model as in the statement and let $\{\delta_i:i\in\mathbb{N}\}$ enumerate the first infinitely many Woodin cardinals of $M$. As usual, $M$ is a model of the Axiom of Choice. 
Since $M$ is a model of Powerset, we have
$$M \models \text{\quotesleft for every $\kappa$, $\kappa^+$ exists.''}$$
A consequence of this is that $M \models \kp$; thus $M$ can define the class $(L(A))^M$ for all $A \in M$.
Assume without loss of generality that $M$ is minimal, in the sense that
$$M \models V = L(V_\lambda).$$
Since $M \models \mathsf{ZC}$, it can form ultrapowers of itself in the usual (non-fine-structural) sense and realize them e.g., as the union of all ultrapowers of initial segments of its cumulative hierarchy. One can carry out the usual proofs of basic facts about the stationary tower $\mathbb{Q}_{<\delta}$, as well as the proof of Woodin's derived model theorem (e.g., the one in \cite{Sta}) within $M$.\footnote{This is where we use the assumption that $\alpha$ is large. For instance, the proof in \cite{Sta} involves building conditions for the stationary tower $\mathbb{Q}_{<\delta}$ starting from elementary substructures of $V_{\lambda+\omega}$. We have not attempted to prove \eqref{LemmaBorelDMMaineq} from optimal hypotheses; it seems plausible that a more careful argument shows e.g., that $L_{\vartheta_1}(\mathbb{R}) \models \ad$ if there is a countable, $\omega_1$-iterable extender model of $\zf$ - Powerset with a limit of Woodin cardinals.} 
Suppose $g^* \subset\coll(\omega,\lambda)$ is $M$-generic. Since $M \models \mathsf{ZC} + \kp$, then 
$$M[g^*] \models \mathsf{ZC} + \kp,$$
(this follows, e.g., from work of Mathias \cite{Mat15}) so that if $\mathbb{R}^*$ is the set of reals of a symmetric collapse at $\lambda$, then
\begin{equation}\label{eqDMinM0}
L(\mathbb{R}^*)^{M}\models \mathsf{Z} + \kp
\end{equation}
(for the proof of this in the case of $L$, see Gostanian \cite[Theorem 1.6]{Go80}.)
Since $\lambda$ is a limit of Woodin cardinals in $M$,
\begin{equation}\label{eqDMinM}
L(\mathbb{R}^*)^{M}\models\ad + \dc.
\end{equation}
Now let $g\subset\coll(\omega,\mathbb{R})$ be $V$-generic. In $V[g]$ let $\{x_i:i\in\mathbb{R}\}$ enumerate $\mathbb{R}$ and find a sequence $M_i$ of extender models such that for each $i\in\mathbb{N}$,
\begin{enumerate}
\item $M_{i+1}$ is the direct limit of an iteration tree on $M_i$ of length $\omega$ by extenders with critical point above the $i$th Woodin cardinal of $M_i$ and length below the $i+1$th Woodin cardinal of $M_{i+1}$;
\item If $j_i:M_0\to M_{i+1}$ is the embedding, then $j_i \in V$ and there is an $M_{i+1}$-generic $g_{i}\subset\coll(\omega,j(\delta_i))$ such that $x_i \in M_{i+1}[g_i]$. 
\end{enumerate}
If so, then, letting $M_\infty$ be the direct limit of $\{M_i:i\in \mathbb{N}\}$, it follows by the restriction on the extenders allowed in the trees that $M_\infty$ is wellfounded. If $j$ is the embedding, then $j(\lambda) = \omega_1^V$. By \eqref{eqDMinM} and the homogeneity of the symmetric collapse, 
$$M \models \text{\quotesleft  $\Vdash$ every set in $L(\dot{\mathbb{R}})$ is determined and $L(\dot{\mathbb{R}})\models\dc$,''}$$
where $\dot{\mathbb{R}}$ is a name for the set of reals of the symmetric collapse. By elementarity, 
\begin{equation}\label{eqDMinM2}
L(\mathbb{R}^*)^{M_\infty}\models\ad+\dc.
\end{equation}
The symmetric collapse can be chosen (in $V[g]$) in such a way that it absorbs the generic $g_i$ for each $i\in\mathbb{N}$, in which case $(\mathbb{R}^*)^{M_\infty} = \mathbb{R}^V$, so $L(\mathbb{R}^*)^{M_\infty}$ is of the form $L_\xi(\mathbb{R})$. 
Suppose $h\subset \coll(\omega,\lambda)$ is $M_\infty$-generic and chosen so that the symmetric collapse induced by $h$ absorbs $g_i$ for each $i\in\mathbb{N}$. By the remark before equation \eqref{eqDMinM0}, we have
$$M_\infty[h] \models \mathsf{ZC}.$$
Since $\coll(\omega,\lambda)$ has the $\lambda^+$-chain condition, 
\begin{equation}\label{eqDMinM1}
M_\infty[h] \models \text{\quotesleft $V_\alpha$ exists.''}
\end{equation}
Now, clearly, we have
$$L_\xi(\mathbb{R})\subset M_\infty[h],
$$ 
so $L_\xi(\mathbb{R})$ satisfies $\ad$ by  \eqref{eqDMinM2}; it satisfies $\mathsf{Z}$ by \eqref{eqDMinM0}, and it satisfies \quotesleft $V_\alpha$ exists'' by \eqref{eqDMinM1}. By minimality, $\vartheta_\alpha<\xi$, which completes the proof.
\endproof

\begin{lemma}\label{LemmaBorelIDM}
Suppose that $L_{\vartheta}(\mathbb{R})\models\ad$. Then, for each $\alpha<\omega_1$, there is an $\omega_1$-iterable extender model of
$$\mathsf{Z} + \text{\quotesleft$V_{\lambda+\alpha}$ exists and $\lambda$ is a limit of Woodin cardinals.''}$$
\end{lemma}
\proof
Suppose that $L_{\vartheta}(\mathbb{R})\models\ad$. Let $\alpha<\omega_1$ be large and find some $\gamma<\vartheta$ such that
\begin{enumerate}
\item $L_\gamma(\mathbb{R})$ is a model of $\mathsf{Z}$ +\quotesleft$V_{\alpha}$ exists,''
\item if $x,y\in\mathbb{R}$ and $x$ is definable in $L_\gamma(\mathbb{R})$ from $y$ and ordinal parameters, then there is $\gamma'<\gamma$ such that $x$ is definable in $L_{\gamma'}(\mathbb{R})$ from $y$ and ordinal parameters.
\end{enumerate}
For example, one could take $\gamma$ to be least such that 
$$V_{\alpha+\omega}^{L_{\vartheta_{\alpha+\omega}}(\mathbb{R})} \subset L_\gamma(\mathbb{R}).$$
Clearly $L_\gamma(\mathbb{R})$ is a model of \quotesleft$V_{\alpha}$ exists.'' Moreover, if $a \in L_{\beta}(\mathbb{R})$, with $\beta<\gamma$, then---by minimality---there is $n\in\mathbb{N}$ such that 
$$V_{\alpha+n}^{L_{\vartheta_{\alpha+\omega}}(\mathbb{R})} \not\subset L_\beta(\mathbb{R}).$$
Thus, there is a stage $\beta'+1 \in (\beta,\gamma)$ at which a new element of $V_{\alpha+n}$ is constructed and, by Proposition \ref{PropositionProjectum}, there is a surjection from $V_{\alpha+n}\cap J_{\beta'}(\mathbb{R})$ to $J_{\beta'}(\mathbb{R})$ in $J_{\beta'+1}(\mathbb{R})$. In particular, there is a surjection from $V_{\alpha+n}\cap J_{\beta'}(\mathbb{R})$ to $a$. Hence, $\mathcal{P}(a)^{L_{\vartheta_{\alpha+\omega}}(\mathbb{R})}\subset V_{\alpha+n+1}^{L_{\vartheta_{\alpha+\omega}}(\mathbb{R})}$ and so 
$$\mathcal{P}(a)^{L_{\vartheta_{\alpha+\omega}}(\mathbb{R})} = \mathcal{P}(a)^{L_\gamma(\mathbb{R})} \in L_\gamma(\mathbb{R}),$$
so we have that $L_\gamma(\mathbb{R})\models \mathsf{Z}$.
Lastly, the definition implies that $\gamma$ is a cardinal (in fact, the largest cardinal) of $L_{\vartheta_{\alpha+\omega}}(\mathbb{R})$ greater than $\Theta^{L_{\vartheta_{\alpha+\omega}}(\mathbb{R})}$, so that (by the L\"owenheim-Skolem theorem),
we have
$$L_{\gamma}(\mathbb{R}) \prec_1 L_{\vartheta_{\alpha+\omega}}(\mathbb{R}).$$
Hence, if $x,y\in\mathbb{R}$ and $x$ is definable in $L_\gamma(\mathbb{R})$ from $y$ and ordinal parameters, then 
$$L_{\gamma}(\mathbb{R})\models \exists \gamma'\, \big(\text{$x$ is definable in $L_{\gamma'}(\mathbb{R})$ from $y$ and ordinal parameters}\big),$$
so $\gamma$ is as desired.

The first thing to verify is that 
$$L_{\gamma}(\mathbb{R})\models\text{\quotesleft Mouse Capturing.''}$$
Recall that \quotesleft Mouse Capturing'' is the statement that for every pair of reals $(x,y)$, $x$ is definable from $y$ and ordinal parameters if, and only if, $x$ belongs to a countably iterable extender model over $y$. In Schindler-Steel \cite[Theorem 3.4.6]{SchSt}, it is shown that for all reals $x,y$, the following are equivalent:
\begin{enumerate}
\item $x$ is definable from $y$ and ordinal parameters in $L_\beta(\mathbb{R})$ for some $\beta<\gamma$,
\item $x$ belongs to an extender model over $y$ which is countably iterable by an iteration strategy in $L_\gamma(\mathbb{R})$.
\end{enumerate}
The equivalence is stated modulo an inductive hypothesis $W^*_\alpha$, but this hypothesis is proved (see the argument in pp. 157-158 of \cite{SchSt}) under the assumption of $\ad$ (in $L(\mathbb{R})$, but the argument goes through in $L_\gamma(\mathbb{R})$). By our choice of $\gamma$, this shows that $L_\gamma(\mathbb{R})$ satisfies Mouse Capturing.

The second thing to verify is that 
$$L_{\gamma}(\mathbb{R})\models \ad^+;$$
this follows from the usual proof for $L(\mathbb{R})$: Ordinal Determinacy follows from the facts that by Kechris-Kleinberg-Moschovakis-Woodin \cite[Theorem 1.1]{KKMW}, $(\delta^2_1)^{L_\gamma(\mathbb{R})}$ is a limit of cardinals with the strong partition property, and that by a theorem due independently to Moschovakis and Woodin (see Larson \cite{La17} for a proof), this implies ${<}(\delta^2_1)^{L_\gamma(\mathbb{R})}$-determinacy. Since
$$L_{(\delta^2_1)^{L_\gamma(\mathbb{R})}}(\mathbb{R})\prec_1 L_\gamma(\mathbb{R}),$$
this implies Ordinal Determinacy in $L_\gamma(\mathbb{R})$.

That every set of reals is $\infty$-Borel in $L_\gamma(\mathbb{R})$ follows from the argument in \cite{La17}. (It uses the fact that $(\Sigma^2_1)^{L_\gamma(\mathbb{R})}$ has the scale property, which is proved by following the argument of \cite{MaSt08} as in Lemma \ref{LemmaScaleVartheta}).

We have checked that $L_\gamma(\mathbb{R})$ satisfies $\ad^+$ and Mouse Capturing. Now, let $\beta<\alpha$ be arbitrary and let $x_\beta$ be a real such that $\beta$ is recursive in $x_\beta$. Define the theory
\[S \equiv \mathsf{Z} + \kp + \ad^+ + \dc + \text{``Mouse Capturing''} + \text{``$V_\beta$ exists.''}\]
Let $\kappa$ be least such that $J_\kappa(\mathbb{R})\models S$; clearly $\kappa<\gamma$ and, in fact, $\kappa<(\delta^2_1)^{L_\gamma(\mathbb{R})}$.
We now run the proof of Steel-Woodin \cite[Theorem 7.2]{StW16} for the theory $S$ within $L_\gamma(\mathbb{R})$ (so that $\kappa$ here takes the role of the ordinal denoted by $\gamma$ therein), except that we consider the pointclass $(\Sigma^2_1(x_\beta))^{J_\kappa(\mathbb{R})}$, rather than $(\Sigma^2_1)^{J_\kappa(\mathbb{R})}$. Similarly, one replaces the $L[E]$-construction on p. 325 of \cite{StW16} by an $L[E](x_\beta)$-construction.
Letting $M$, $N$, and $\xi$ be as in \cite{StW16}, and $M_0$ be set of all elements of $M|\xi$ which are definable \emph{from $x_\beta$} in $M|\xi$, the argument shows that $M_0$ has the form $J_{\zeta}(N_0)$, for some $\zeta$ and some extender model $N_0$ over $x_\beta$ and, moreover, the following hold:
\begin{enumerate}
\item $N_0$ satisfies ``there are infinitely many Woodin cardinals'';
\item no initial segment of $M_0$ projects to $N_0$;
\item $M_0$ has an $(\omega_1,\omega_1)$-iteration strategy $\Sigma_0$;
\item the derived model of $M_0$ at its limit of Woodin cardinals satisfies $S$.
\end{enumerate}
Now, let $\lambda$ be the limit of Woodin cardinals of $M_0$. Then, an ordinal $\eta>\lambda$ is a cardinal of $M_0$ if, and only if, it is a cardinal of a generic extension of $M_0$ by $\coll(\omega,\lambda)$ if, and only if, it is a cardinal of the derived model of $M_0$. Since this derived model satisfies $S$, it follows that $M_0$ has $\beta$-many cardinals above $\lambda$ and, since $M_0$ is a model of the Generalized Continuum Hypothesis, we have $M_0 \models$ ``$V_\beta$ exists.''
Since $\beta<\alpha$ was arbitrary and $\alpha$ was arbitrary, the result follows.
\endproof

Theorem \ref{mainborelwoodins} is immediate from Theorem \ref{mainborel},  Lemma \ref{LemmaBorelIDM} and Lemma \ref{LemmaBorelDM}.

\bibliographystyle{abstract}
\bibliography{References}

\end{document}